\documentclass[12pt]{article}

\usepackage[margin=1in]{geometry}  
\usepackage{graphicx}              
\usepackage{amsmath}               
\usepackage{amsfonts}              
\usepackage{amsthm}                
\usepackage{amssymb}

\begin{document}

\nocite{*}

\title{Intersection lattices of cubic fourfolds.}

\author{Evgeny Mayanskiy}

\maketitle

\begin{abstract}
  We study the lattices of algebraic and transcendental cycles of cubic fourfolds.
\end{abstract}

\section{Introduction}

\subsection{Motivation.}

Our objective is to understand which lattices occur as sublattices of algebraic and transcendental cycles of cubic fourfolds. The analogous problem for $K3$ surfaces has been solved as an application of Nikulin's results \cite{Nikulin} and the nice properties (surjectivity) of the period map for $K3$ surfaces \cite{BPV}. Its statement and solution can be found, for example, in \cite{Morrison}.\\

Such results may be applied to the study of automorphism groups of these varieties. Namely, we know from the Torelli theorem for $K3$ surfaces \cite{BPV} that the automorphism group of a $K3$ surface is isomorphic to the group of effective Hodge isometries of its middle integral cohomology lattice. As for cubic fourfolds, Torelli theorem due to Voisin \cite{Voisin} says that every Hodge isometry of the middle integral cohomology lattice of a cubic fourfold $X$, preserving the square of the hyperplane section class, comes from an automorphism of $X$. \\

The question about automorphism groups of cubic fourfolds will be studied in our future work. \\

Another possible application is related to the rationality question for cubic fourfolds. So far no example of a non-rational cubic fourfold was found. New examples of rational cubic fourfolds would be of interest too. The relevance of the study of intersection lattices of cubic fourfolds to such problems follows from results of several people, in particular, from the work by Hassett \cite{Hassett}, \cite{Hassett2}. We hope that our study will facilitate further advance in this direction.  

\subsection{Main results.}

For a cubic fourfold $X$, we will denote by $A(X)$ the lattice generated by algebraic cycles (i.e. $A(X)=H^4(X,\mathbb Z)\cap H^{2,2}(X)$), by $T(X)$ its lattice of transcendental cycles (i.e. $T(X)=(A(X))^{\bot}_{H^4(X,\mathbb Z)}$) and by $A_0(X)$ the lattice of primitive algebraic cycles (i.e. the orthogonal complement in $A(X)$ to the class of a hyperplane section). These lattices are primitive sublattices of the middle integral cohomology lattice
$$
H^4(X, \mathbb Z)\cong L \mathrel{\mathop:}= (1)^{\oplus 21} \oplus (-1)^{\oplus 2},
$$
and its primitive piece $H^4(X,\mathbb Z)_0$ which is isomorphic to
$$
L_0 \mathrel{\mathop:}= \begin{pmatrix} 2 & 1\\ 1 & 2 \end{pmatrix} \oplus U^{\oplus 2} \oplus E_8^{\oplus 2}.
$$
Given a lattice $S$ with quadratic form $Q$, we will call $x\in S$ a {\it root}, if $Q(x,x)=2$. We will call $y\in S$ a {\it long root}, if $Q(y,y)=6$ and $Q(y,S)\subset 3\mathbb Z$.\\ 

For the other notation and terminology we refer the reader to Sect. 2.\\

The rest of this section contains the main results of the paper. In all the statements a positive integer $\rho$ coincides with the rank of the lattice $A(X)$ of algebraic cycles on the corresponding cubic fourfold $X$. In particular, when $\rho (X)=1$, $A(X)=\mathbb Z h^2$ (a rank $1$ lattice with matrix $< 3 >$), $A_0(X)=0$ and $T(X)=H^4(X,\mathbb Z)_0\cong L_0$.\\

{\bf Corollary 5.1} {\it Let $2\leq \rho \leq 11$, and $A_0$ be an even lattice of signature $(\rho-1,0)$ that does not contain roots. Then there is a cubic fourfold $X$ and an isomorphism of lattices $A_0(X) \cong A_0$.}\\ 

{\bf Corollary 5.2} {\it Let $13\leq \rho \leq 21$, and $T$ be an even lattice of signature $(21 - \rho,2)$. Then there is a primitive embedding of $T$ into $L_0$. Suppose one can choose it in such a way that the orthogonal complement of $T$ in $L_0$ contains neither roots, nor long roots of $L_0$. Then there is a cubic fourfold $X$ and an isomorphism of lattices $T(X)  \cong T$.}\\ 

For a finite abelian group $B$ we define $l(B)$ to be the minimal number of generators of $B$ and $B_p$ to be its p-component (in the decomposition of $B$ into a direct sum of finite abelian p-groups $B={\bigoplus}_p B_p$). In the next corollary we use notation from the classification of p-adic lattices \cite{Nikulin}, which is recalled in section 2.3 below. In particular, we use invariants $v_p\in \{ 0,1 \}$ and ${\theta}_p \in {\mathbb Z}_{p}^{*}/({\mathbb Z}_{p}^{*})^2$.\\ 

{\bf Corollary 5.3} {\it Let $S$ be an even lattice without roots, $sign(S)\leq sign(L_0)$. Let $q \colon A_S \rightarrow {\mathbb Q}/2\mathbb Z$ be the corresponding discriminant form. Then there exists a cubic fourfold $X$ and an isomorphism of lattices $S\cong A_0(X)$, if and only if one of the following conditions holds:
\begin{enumerate}
				\item[A1] $l(A_S)=l((A_S)_3)\leq 20-rk(S)$ 
				\item[A2] $l(A_S)=l((A_S)_3)=21-rk(S)$ and $\left( \frac{(-1)^{rk(S)}\cdot ({\theta}_3)^{v_3}}{3} \right)=1$
				\item[A3] {$22-rk(S)=l(A_S)>l((A_S)_3)$ and the following conditions hold: 
						\begin{itemize}
						\item If for a prime $p\neq 2,3$, $l(A_S)=l((A_S)_p)$, then $\left( \frac{(-1)^{rk(S)}\cdot 3\cdot ({\theta}_p)^{v_p}}{p} \right)=1$.	
						\item If $l(A_S)=l((A_S)_3)+1$, then $\left( \frac{(-1)^{rk(S)+1}\cdot ({\theta}_3)^{v_3}}{3} \right)=1$.	
						\item If $l(A_S)=l((A_S)_2)$ and $q_2\neq q_{\theta}^{(2)}(2)\oplus q_2^{'}$, then $v_2=1$ (in this case also $rk(S)>l(A_S)$).\\
						\end{itemize} }
\item[B] There exists an element $h\in A_S$ of order $3$ such that $q_S(h)=\frac{2}{3}\in \mathbb Q/2 \mathbb Z$ and does not exist $\delta \in S$ with the property that $Q_S(\delta, \delta)=6$ and $\sigma (\delta)\in 3 \mathbb Z$ for any $\sigma \in S^{*}$ such that $Q_{S^{*}}(\sigma, h)\in \mathbb Z$, and one of the following conditions holds: 
				\begin{enumerate}
				\item[B1] $l(A_S)=l((A_S)_3)\leq 22-rk(S)$ 
				\item[B2] $l(A_S)=l((A_S)_3)=23-rk(S)$ and the following conditions hold: 
						\begin{itemize}
						\item If for a prime $p\neq 2,3$, $l(A_S)=1+l((A_S)_p)$, then $\left( \frac{(-1)^{rk(S)}\cdot 3\cdot ({\theta}_p)^{v_p}}{p} \right)=1$.	
						\item If $l(A_S)=l((A_S)_3)$, then $\left( \frac{(-1)^{rk(S)}\cdot ({\theta}_3)^{v_3}}{3} \right)=1$.	
						\item If $l(A_S)=1+l((A_S)_2)$ and $q_2\neq q_{\theta}^{(2)}(2)\oplus q_2^{'}$, then $v_2=1$ (in this case also $rk(S)>l(A_S)$).	
						\end{itemize} 
				\item[B3] $22-rk(S)=l(A_S)>l((A_S)_3)$ and the following conditions hold: 
						\begin{itemize}
						\item If for a prime $p\neq 2,3$, $l(A_S)=l((A_S)_p)$, then $\left( \frac{(-1)^{rk(S)}\cdot 3\cdot ({\theta}_p)^{v_p}}{p} \right)=1$.	
						\item If $l(A_S)=l((A_S)_2)$ and $q_2\neq q_{\theta}^{(2)}(2)\oplus q_2^{'}$, then $v_2=1$ (in this case also $rk(S)>l(A_S)$).\\
						\end{itemize} 
				\end{enumerate}
\end{enumerate}}

{\bf Theorem 6.1} {\it Let $A$ be an odd lattice of signature $sign(A)=(\rho,0)$, where $2\leq \rho \leq 21$. Then $A\cong A(X)$ for some cubic fourfold $X$ if and only if there exists $a\in A$, such that if we denote $A_0=(\mathbb Z a)_{A}^{\bot}$, then the following conditions are satisfied:
\begin{enumerate}
\item $Q_A(a,a)=3$.
\item $A_0$ is an even lattice.
\item There is {\it no} ${\delta}_0\in A_0$ such that $Q_A({\delta}_0,{\delta}_0)=2$.
\item There is {\it no} ${\delta}\in A_0$ such that $Q_A({\delta},{\delta})=6$ and for any $\alpha \in A^{*}$ with ${\alpha}(a)=0$ we have ${\alpha}(\delta)\in 3 \mathbb Z$.
\item For any $b\in A$ the integer $Q_A(a,b)^2-Q_A(b,b)$ is even.
\item There exists an even lattice $K$ of signature $sign(K)=sign(L)-sign(A)=(21-\rho,2)$, such that $A_A=A_K$ and $q_K \colon A_A=A^{*}/A \rightarrow \mathbb Q/ 2\mathbb Z$, $\alpha \mapsto \left[ ({\alpha}(a))^2-Q_{A^{*}}(\alpha,\alpha) \right] +2\mathbb Z$.
\end{enumerate}
}

{\bf Theorem 5.4} {\it Let $2\leq \rho \leq 21$, $T$ be an even lattice with signature $sign(T)=(21-\rho,2)$ and discriminant-form $q_T\colon A_T \rightarrow \mathbb Q/ 2\mathbb Z$. Then there exists a cubic fourfold $X$ and an isomorphism of lattices $T(X)\cong T$, if and only if one of the following conditions holds:
\begin{enumerate}
\item[A] There exists $\tau\in A_T$ an element of order $3$ such that $q_T(\tau)=\frac{2}{3}+2\mathbb Z\in \mathbb Q / 2 \mathbb Z$, and if we denote $G=(\tau)^{\bot}_{A_T}\subset A_T$, then there exists an even lattice $K$ of signature $(\rho-1,0)$ such that $A_K=G\subset A_T$, $q_K=-q_T {\mid}_{G}$, and $K$ has no roots.
\item[B] There exists a primitive embedding of lattices $U^{23-\rho} \hookrightarrow T(-1)\oplus L_0$ such that the orthogonal complement $K$ of its image has no roots, and for any $\Delta \in (T(-1)\oplus {L_0}^{*})  \backslash (T(-1)\oplus {L_0})$ such that ${\Delta}(U^{23-\rho})=0$, ${\Delta}(3{\Delta})=\frac{1}{3}Q_{T(-1)\oplus {L_0}}(3{\Delta},3{\Delta})\neq 2$.
\end{enumerate}}

\subsection{Acknowledgements.}
Our work uses essentially the results of Nikulin \cite{Nikulin}, Laza \cite{Laza} and Looijenga \cite{Looijenga}. We also benefited greatly from reading the paper by Hassett \cite{Hassett}. We understood a possible general approach to this problem after reading the paper by Morrison \cite{Morrison}, where the solution of the analogous problem for $K3$ surfaces may be found. We also appreciate the help and advice we have been obtaining during the course of our PhD study from our advisor Yuri Zarhin and the support from the Mathematical department of the Pennsylvania State University. In addition, we thank Yuri Zarhin for assigning this problem to us.\\

\section{Notation.}

Throughout the paper all varieties are complex and smooth.  

\subsection{General notation.}

Let $E_8$ and $U$ be even unimodular lattices of signatures $(8,0)$ and $(1,1)$ respectively.\\
We will also use the following lattices: 
$$
L=(1)^{\oplus 21} \oplus (-1)^{\oplus 2},
$$

$$
L_0= \begin{pmatrix} 2 & 1\\ 1 & 2 \end{pmatrix} \oplus U^{\oplus 2} \oplus E_8^{\oplus 2}.
$$
Let $Q$ denote the bilinear form of $L$ and $L_0$.\\

It is known that the middle integral cohomology lattice and the primitive middle cohomology lattice of a cubic fourfold $X$ are isomorphic to $L$ and $L_0$ respectively \cite{Hassett}: 
$$
H^4(X,\mathbb Z)\cong L {\mbox \qquad and \qquad } H^4(X,\mathbb Z)_0 \cong L_0
$$
For a cubic fourfold $X$ we denote by $A_0(X)$ the orthogonal complement to the self-intersection of the hyperplane section class in its saturated lattice of algebraic cycles:
$$
A_0(X)=H^{4}(X, \mathbb Z)_0 \cap H^{2,2}(X)  
$$
and by $T(X)$ the transcendental lattice of $X$, i.e. the orthogonal complement of $A_0(X)$ in $H^4(X,\mathbb Z)_0$:
$$
T(X) \cong A_0(X)^{\bot}
$$
Let's choose and fix an element $h^2 \in L$ of self-intersection $Q(h^2,h^2)=3$, and a primitive embedding $L_0 \hookrightarrow L$, identifying $L_0$ with the orthogonal complement of $h^2$.\\

{\bf Remark.} It follows form the theory of discriminant-forms \cite{Nikulin} (namely, from the odd analog of Theorem 1.14.4 and the discussion in section 16) that the primitive embedding $L_0 \hookrightarrow L$ is unique upto an isomorphism, and the orthogonal complement of $L_0$ in $L$ is always generated by an element $h^2\in L$ with self-intersection 3.\\

Symbol $\rho$ denotes an integer with $1 \leq \rho \leq 21$. It will denote the rank of the (saturated) sublattice of algebraic cycles of a cubic fourfold $X$:
$H^{4}(X, \mathbb Z) \cap H^{2,2}(X)$. The range $1 \leq \rho \leq 21$ is explained and justified by the result of Zarhin \cite{Zarhin}, which says that for any such $\rho$ there exists a cubic fourfold $X$, such that $rk (H^{4}(X, \mathbb Z) \cap H^{2,2}(X))=\rho$.\\

\subsection{Period domain and period map.}

Following \cite{Hassett}, we recall the definition of the period domain and the period map for cubic fourfolds.\\

Let $\mathcal C$ be the coarse moduli space of cubic fourfolds and $\mathcal D'$ be the local period domain, i.e. one of the two connected components of the analytic open subset of the quadric 
$$
\{  x\in \mathbb P (L_0 \otimes_{\mathbb Z} \mathbb C) | Q(x,x)=0, Q(x,\bar x)<0 \} 
$$
Then the global period domain for cubic fourfolds is denoted by $\mathcal D$ and is defined as $\mathcal D=\Gamma^{+} \backslash \mathcal D' $, where $\Gamma^{+}$ is the group of automorphisms of lattice $L$, preserving $h^2$ and the chosen connected component of the quadric above.\\ 

{\bf Theorem [Voisin, Hassett]} {\it The period map for cubic fourfolds $\tau \colon \mathcal C \longrightarrow \mathcal D$ is an open embedding of quasiprojective varieties of dimension 20.}\\

\subsection{Lattice theoretic notation.}

We will use the same lattice theoretic notation as in \cite{Nikulin} and \cite{Morrison}.\\

In particular, for a finite group $G$, we denote by $l(G)$ the minimal possible number of its generators, and for an even lattice $S$ its discriminant form $q_S \colon A_S \rightarrow \mathbb Q / 2 \mathbb Z$, where $A_S=S^{*}/S$. Notice that $l(A_S) \leq rk(S)$. If $A_S={\oplus}_p \oplus_{k\geq 1}(\mathbb Z/p^k \mathbb Z)^{c_k(p)}$, then $l(A_S)={max}_p \{ \sum_{k\geq 1} c_k(p)\}$.\\ 

For the reader's convenience let us recall the construction of the finite symmetric nondegenerate bilinear form $b_S\colon A_S \otimes A_S \rightarrow \mathbb Q / \mathbb Z$ and the discriminant-form $q_S \colon A_S \rightarrow \mathbb Q / 2 \mathbb Z$, associated to a lattice $S$ (in the definition of $q_S$ we assume that $S$ is even). For details, see \cite{Nikulin}.\\

The bilinear form $Q_S \colon S\otimes S \rightarrow \mathbb Z$ on the lattice $S$ extends (uniquely) to a bilinear form $Q_{S^{*}} \colon S^{*}\otimes S^{*} \rightarrow \mathbb Q$ on the free abelian group $S^{*}$. (Note that there exists a canonical embedding of abelian groups $S\hookrightarrow S^{*}$.) Then we define:
$$
b_S \colon (S^{*}/S)\otimes (S^{*}/S)\rightarrow \mathbb Q / \mathbb Z \mbox{,} \quad  \quad (x,y)\mapsto Q_{S^{*}}(x,y)+\mathbb Z
$$
and (if $S$ is even)
$$
q_S \colon (S^{*}/S)\rightarrow \mathbb Q /2 \mathbb Z \mbox{,} \quad  \quad x \mapsto Q_{S^{*}}(x,x)+2 \mathbb Z.
$$

For a lattice $S$ we call $x\in S$ a {\it root}, if $Q_S(x,x)=2$, and a {\it long root}, if $Q_S(x,x)=6$ and $Q_S(x,y)\in 3 \mathbb Z$ for all $y \in S$.\\

A sublattice $S \subset N$ is called {\it saturated}, if $S=N\cap (S{\otimes}_{\mathbb Z}{\mathbb Q})\subset N{\otimes}_{\mathbb Z}{\mathbb Q}$. For any sublattice $S \subset N$ there exists a minimal saturated sublattice $S^{\wedge}$ of $N$, containing $S$. Sublattice $S^{\wedge}\subset N$ is called the {\it saturation} of $S$ in $N$.\\

In our classification we use local invariants of even lattices. Let us recall classification of even p-adic lattices following \cite{Nikulin}.\\

Let $K_{\theta}^{(p)}(p^k)$ denote the rank $1$ p-adic lattice with matrix $< {\theta}p^k >$, where $k\geq 0$ and ${\theta}\in {\mathbb Z}_{p}^{*}/({\mathbb Z}_{p}^{*})^2$. Let $U^{(2)}(2^k)$ and $V^{(2)}(2^k)$ ($k\geq 0$) be rank $2$ 2-adic lattices with matrices
$\begin{pmatrix}
0 & 2^k \\
2^k & 0 \end{pmatrix}$ and 
$\begin{pmatrix}
2^{k+1} & 2^k \\
2^k & 2^{k+1} \end{pmatrix}$ 
respectively. Then we will denote by $q_{\theta}^{(p)}(p^k)$, $u_{+}^{(2)}(2^k)$ and $v_{+}^{(2)}(2^k)$ their discriminant-forms.\\

{\bf Theorem.} (\cite{Nikulin}, Theorem 1.9.1) {\it Let $A_p$ be a finite abelian p-group and $q_p\colon A_p \rightarrow {\mathbb Q}_p /2{\mathbb Z}_p$ a finite non-degenerate quadratic form. Then there exists a unique even p-adic lattice $K(q_p)$ of rank $l(A_p)$, whose discriminant-form is isomorphic to $q_p$, unless $p=2$ and $q_2=q_{\theta}^{(2)}(2)\oplus q_2^{'}$. In the latter case, there are exactly two even 2-adic lattices $K_{{\alpha}_1}(q_2)$ and $K_{{\alpha}_2}(q_2)$, whose  discriminant-forms are isomorphic to $q_2$.}\\

{\bf Theorem.} (\cite{Nikulin}, Corollary 1.9.3) {\it Let $K_p$ be an even p-adic lattice with discriminant-form $q_p$. Then $K_p$ can be uniquely represented in the following form:
\begin{enumerate}
\item[(a)] If $p\neq 2$, then $K_p\cong K_1^{(p)}(1)^{t_p-v_p}\oplus K_{{\theta}_p}^{(p)}(1)^{v_p}\oplus K(q_p)$, where $\left( \frac{{\theta}_p}{p} \right)=-1$, $0\leq v_p\leq 1$, $t_p\geq v_p$.
\item[(b)] If $p=2$, $q_2\neq q_{\theta}^{(2)}(2)\oplus q_2^{'}$, then $K_2\cong U^{(2)}(1)^{t_2-v_2}\oplus V^{(2)}(1)^{v_2}\oplus K(q_2)$, where $0\leq v_2\leq 1$, $t_2\geq v_2$.
\item[(c)] If $p=2$, $q_2= q_{\theta}^{(2)}(2)\oplus q_2^{'}$, then $K_2\cong U^{(2)}(1)^{t_2}\oplus K(q_2)$, where $t_2\geq 0$.\\
\end{enumerate}}

Invariants ${\theta}_p \in {{\mathbb Z}_{p}^{*}}/{({\mathbb Z}_{p}^{*})^2}$ and $v_p\in \{ 0, 1 \}$ of the p-adic lattice $S_p$ obtained by localization of an even lattice $S$ will be essential for our classification. Given an even lattice $S$, consider the corresponding finite abelian group $A=A_S=S^{*}/S$ and the discriminant-form $q=q_S\colon A\rightarrow \mathbb Q/2\mathbb Z$. If $A={\oplus}_p A_p$, $q={\oplus}_p q_p$ is the decomposition into a direct sum of p-groups, then the p-adic lattice corresponding to $S$ is the lattice $S_p=S{\otimes}_{\mathbb Z} {\mathbb Z}_p$ with discriminant-form $q_p \colon A_p \rightarrow {\mathbb Q}_p /2 {\mathbb Z}_p$.  \\

{\bf Remark 2.1 } Let $S$ be an even lattice with discriminant-form $q\colon A\rightarrow \mathbb Q/2 \mathbb Z$, where $A=A_S=S^{*}/S$. Then by the Theorem above (Corollary 1.9.3 of \cite{Nikulin}) for any odd prime $p$ we have $discr(S)=discr(S_p)={{\theta}_p}^{v_p}\cdot discr(K(q_p))$, and so $\frac{(-1)^{t_{+}}\cdot \mid A \mid}{discr(K(q_p))}=(-1)^{rk(S)}\cdot {{\theta}_p}^{v_p}$. (Notice that $discr(S)=(-1)^{t_{-}}\cdot \mid A \mid$ and $t_p=rk(S)-l(A_p)$.)\\

If $q_2\neq q_{\theta}^{(2)}(2)\oplus q_2^{'}$, then $t_2=(rk(S)-l(A_2))/2$ and $discr(S)=discr(S_2)=(-1)^{\frac{rk(S)-l(A_2)}{2}}\cdot 3^{v_2}\cdot discr(K(q_2))$, and so $\frac{\mid A \mid}{discr(K(q_2))}=(-1)^{t_{-}}\cdot (-1)^{(rk(S)-l(A_2))/2} \cdot 3^{v_2}$. Let us note that when $v_2=1$, we should have $l(A_2)< rk(S)$.\\ 

\section{Existence of cubic fourfolds for a given primitively embedded lattice.}

We will use the following result of Laza and Looijenga, which we state in the form of Theorem 1.1 of paper \cite{Laza}.\\ 

{\bf Theorem 3.1 [Laza, Looijenga]} {\it For a primitive sublattice $M \hookrightarrow L$ of rank $2$, containing $h^2$, and such that $det M=2$ or $det M = 6$, consider the hyperplane $D_M \subset L_0 \otimes_{\mathbb Z} \mathbb C$, defined as 
$$
D_M=\{ x \in \mathbb P (L_0 \otimes_{\mathbb Z} \mathbb C) | Q(x, M)=0 \} 
$$
Then the point $p \in \mathcal D$ does not lie in the image of the period map $\tau$ if and only if there exists $M$ as above such that $p \in D_M$.}\\


Now we can state our existence result.\\

{\bf Corollary 3.2} (Compare with Proposition 2.15 in \cite{Laza}) {\it Let $2 \leq \rho \leq 20$. Let $A_0 \subset L_0$ be a primitive sublattice of signature $(\rho -1,0)$.\\
Then the existence of a cubic fourfold $X$ and an isomorphism of lattices $H^4(X, \mathbb Z) \cong L$, identifying the self-intersection of the hyperplane section class with $h^2$, and $A_0(X)$ with $A_0$, is equivalent to the condition that $A_0$ contains neither roots, nor long roots of $L_0$.}\\

The proof (whose general scheme follows the one of \cite{Morrison}, Corollary 1.9, for $K3$ surfaces) consists of the following observations.\\

{\bf Lemma 3.3} {\it The following statements are equivalent: 
\begin{enumerate}
\item[(a)] $A_0$ contains a root or a long root of $L_0$.
\item[(b)] There exists $m_0\in A_0$ such that $\mathbb Z m_0 \subset A_0$ is saturated and either $det((\mathbb Z m_0 + \mathbb Z h^2)^{\wedge})=2$, or $det((\mathbb Z m_0 + \mathbb Z h^2)^{\wedge})=6$.
\item[(c)] There exists a primitive sublattice $M \subset L$ as in Theorem 3.1 above, such that $(A_0 \otimes_{\mathbb Z} \mathbb C)^{\bot} \subset D_M$ in $L_0 \otimes_{\mathbb Z} \mathbb C$. 
\end{enumerate} 
} 

{\it Proof.} The inclusion in (c) is equivalent to $D_M^{\bot} \subset ((A_0 \otimes_{\mathbb Z} \mathbb C)^{\bot})^{\bot}$ $ = A_0 \otimes_{\mathbb Z} \mathbb C$ in $L_0  \otimes_{\mathbb Z} \mathbb C$.\\

Since $D_M^{\bot}=(M \otimes_{\mathbb Z} \mathbb C)\cap (L_0 \otimes_{\mathbb Z} \mathbb C)$ $=(M \cap L_0) \otimes_{\mathbb Z} \mathbb C = (\mathbb Z h^2)_M^{\bot} \otimes_{\mathbb Z} \mathbb C$ in $L \otimes_{\mathbb Z} \mathbb C$, we conclude that (b) is equivalent to the existence of the primitive sublattice $M \subset L$ as above, such that $(\mathbb Z h^2)_M^{\bot} \subset A_0$.\\ 

Consider a lattice $M_1$ of the form $M_1=(\mathbb Z m_0 + \mathbb Z h^2)^{\wedge}\subset L$, where $m_0 \in A_0$. Then $M_1$ is a primitive sublattice of $L$ of rank $2$, containing $h^2$, such that $(\mathbb Z h^2)_{M_1}^{\bot} \subset A_0$. Any lattice with these properties has such a form for some $m_0 \in A_0$.\\

Hence condition (c) is equivalent to the existence of $m_0 \in A_0$ such that $det((\mathbb Z m_0 + \mathbb Z h^2)^{\wedge})=2$ or $det((\mathbb Z m_0 + \mathbb Z h^2)^{\wedge})=6$. Without loss of generality, we may assume that $\mathbb Z m_0 \subset L$ is saturated. This proves equivalence of (b) and (c).\\

Recall the formula 
$$
\left| \frac{L_1}{M_1} \right| = \sqrt{\left| \frac{det(M_1)}{det(L_1)} \right|} \eqno{(0)}
$$
where $L_1$ is a lattice and $M_1 \subset L_1$ is a sublattice of rank $rk M_1=rk L_1$.\\ 

We get

$$
3 Q(m_0,m_0)=det(\mathbb Z m_0 + \mathbb Z h^2)=det((\mathbb Z m_0 + \mathbb Z h^2)^{\wedge}) \cdot \left|\frac{(\mathbb Z m_0 + \mathbb Z h^2)^{\wedge}}{\mathbb Z m_0 + \mathbb Z h^2}\right|^2 \eqno{(1)}
$$

{\bf Claim 1} For any $\xi \in (\mathbb Z m_0 + \mathbb Z h^2)^{\wedge}$, $3 \xi \in \mathbb Z m_0 + \mathbb Z h^2$.\\

{\it Proof of Claim 1:} Indeed, let $\xi = p m_0 + q h^2$, where $p,q \in \mathbb Q$. Then $3q = Q(\xi, h^2)\in \mathbb Z$. Hence $3p m_0=3 \xi - 3qh^2\in L$, and so $3p \in \mathbb Z$, because $\mathbb Z m_0 \subset L$ is saturated. Hence $3 \xi=(3p)m_0+(3q)h^2\in \mathbb Z m_0 + \mathbb Z h^2$. {\it Claim 1 is proven.}\\ 

{\bf Claim 2} Unless $\mathbb Z m_0 + \mathbb Z h^2$ is already saturated in $L$, $\frac{(\mathbb Z m_0 + \mathbb Z h^2)^{\wedge}}{\mathbb Z m_0 + \mathbb Z h^2}\cong \mathbb Z / 3 \mathbb Z$. \\

{\it Proof of Claim 2:} Indeed, since $\mathbb Z h^2 \subset (\mathbb Z m_0 + \mathbb Z h^2)^{\wedge}$ is a primitive sublattice, we can complete $h^2$ to a $\mathbb Z$-basis $\{ h^2, T \}$ of a $\mathbb Z$-module $(\mathbb Z m_0 + \mathbb Z h^2)^{\wedge}$, i.e. $(\mathbb Z m_0 + \mathbb Z h^2)^{\wedge} = \mathbb Z T + \mathbb Z h^2$.\\

Let $d=det((\mathbb Z m_0 + \mathbb Z h^2)^{\wedge})$ $=det \begin{pmatrix} Q(T,T) & Q(T, h^2)\\ Q(T, h^2) & 3 \end{pmatrix} =$
$$
=3 Q(T,T)- (Q(T, h^2))^2 \eqno{(2)}
$$

Let $m_0=a T + b h^2$, where $a$ and $b$ are integers.\\ 

Then $0=Q(m_0, h^2)=3b+a \cdot Q(T, h^2)$.\\

So, $Q(T, h^2)=-\frac{3b}{a}$, and $Q(T,T)=\frac{d}{3}+\frac{3b^2}{a^2}$.\\

$Q(m_0,m_0)=a^2 Q(T,T)+2ab Q(T,h^2)+3b^2=a^2\cdot (\frac{d}{3}+\frac{3b^2}{a^2})+2ab (-\frac{3b}{a})+3b^2=\frac{a^2d}{3}$.\\

Hence from formula (1) we get $a=\pm \left| \frac{(\mathbb Z m_0 + \mathbb Z h^2)^{\wedge}}{\mathbb Z m_0 + \mathbb Z h^2} \right|$.\\

Now $T=\frac{1}{a}m_0-\frac{b}{a}h^2 \in L$.\\

By the Claim 1, $3T=\frac{3}{a}m_0-\frac{3b}{a}h^2 \in \mathbb Z m_0 + \mathbb Z h^2$. Hence $a=\pm 1$ or $a=\pm 3$.\\

So, we conclude that either $\mathbb Z m_0 + \mathbb Z h^2$ is already satuarted in $L$, or $\frac{(\mathbb Z m_0 + \mathbb Z h^2)^{\wedge}}{\mathbb Z m_0 + \mathbb Z h^2}\cong \mathbb Z / 3 \mathbb Z$.\\

{\it Claim 2 is proven.}\\ 

It follows from Claim 2 and equality (1) above that either

$$
Q(m_0,m_0)=\frac{det(\mathbb Z m_0 + \mathbb Z h^2)}{3} \eqno{(3)}
$$

and $\mathbb Z m_0 + \mathbb Z h^2$ is saturated in $L$, or

$$
Q(m_0,m_0)=3\cdot det((\mathbb Z m_0 + \mathbb Z h^2)^{\wedge}) \eqno{(4)}
$$

and $\frac{(\mathbb Z m_0 + \mathbb Z h^2)^{\wedge}}{\mathbb Z m_0 + \mathbb Z h^2}\cong \mathbb Z / 3 \mathbb Z$.\\ 

{\bf Claim 3} If $det((\mathbb Z m_0 + \mathbb Z h^2)^{\wedge})=6$, then $\mathbb Z m_0 + \mathbb Z h^2$ is saturated in $L$.\\

{\it Proof of Claim 3:} We argue by contradiction. Suppose $\mathbb Z m_0 + \mathbb Z h^2$ is not saturated in $L$. Then by Claim 2 $\frac{(\mathbb Z m_0 + \mathbb Z h^2)^{\wedge}}{\mathbb Z m_0 + \mathbb Z h^2}\cong \mathbb Z / 3 \mathbb Z$, and we are in Case (4), i.e. $Q(m_0,m_0)=18$.\\

Then there is an element $T\in L$ of the form $T=\frac{1}{3}m_0-\frac{b}{3}h^2$. If $3 | b$, then $\frac{1}{3}m_0 \in L$, which contradicts to the assumption that $\mathbb Z m_0 \subset L$ is saturated. Otherwise, an element of the form $T'=\frac{1}{3}(m_0\pm h^2)$ lies in $L$. Hence $Q(T',T')=\frac{1}{9}(Q(m_0,m_0)+Q(h^2,h^2))=\frac{21}{9}$ should be an integer. This is a contradiction. 
 
{\it Claim 3 is proven.}\\ 

{\bf Claim 4} If $m_0 \in A_0$ is a long root of $L_0$, then $\mathbb Z m_0 + \mathbb Z h^2$ is not saturated in $L$.\\

{\it Proof of Claim 4:} We argue by contradiction. Suppose $\mathbb Z m_0 + \mathbb Z h^2$ is saturated in $L$.\\

Since $m_0 \in A_0$ is a long root, $m_0 \mathbb Z \subset L$ is saturated.\\

Let $y \in L$ be a generator of the group $\frac{L}{\mathbb Z + (\mathbb Z h^2)_L^{\bot}} \cong \mathbb Z / 3 \mathbb Z$. Consider an element $\alpha \in Hom_{\mathbb Z}(L, \frac{1}{3}\mathbb Z)$, defined as $\alpha (x)=\frac{1}{3}\cdot Q(am_0+bh^2,x)$ for some pair of integer numbers $a\in \mathbb Z$, $b \in \mathbb Z$ such that at least one of them is not divisible by $3$.\\ 

Then $\alpha \in L^{*}=Hom_{\mathbb Z}(L, \mathbb Z) \subset Hom_{\mathbb Z}(L, \frac{1}{3}\mathbb Z)$ if and only if $\alpha (y)\in \mathbb Z$, i.e. $a$ and $b$ satisfy  relation
$$
a\cdot Q(m_0,y)+b\cdot Q(h^2,y) \equiv 0 (mod 3)
$$

Such $a\in \mathbb Z$ and $b \in \mathbb Z$ always exist, and so the corresponding $\alpha \in Hom_{\mathbb Z}(L, \frac{1}{3}\mathbb Z)$ in fact lies in $L^{*}=Hom_{\mathbb Z}(L, \mathbb Z)$. Since $L$ is unimodular, $L \cong L^{*}$, and so $\frac{1}{3}\cdot (am_0+bh^2)\in (\mathbb Z m_0 + \mathbb Z h^2)^{\wedge} \backslash (\mathbb Z m_0 + \mathbb Z h^2)$. This contradicts to our assumption.\\

{\it Claim 4 is proven.}\\ 

Now let us prove that (b) implies (a).\\

Let $det((\mathbb Z m_0 + \mathbb Z h^2)^{\wedge})=2$.\\

Then $\frac{(\mathbb Z m_0 + \mathbb Z h^2)^{\wedge}}{\mathbb Z m_0 + \mathbb Z h^2}\cong \mathbb Z / 3 \mathbb Z$, and $Q(m_0,m_0)=6$ by formulas (3) and (4).\\

Also, in this case there exists $T=\frac{1}{3}m_0-\frac{b}{3}h^2\in L$, and so $\forall y \in L_0$  $Q(m_0,y)=Q(3T+bh^2,y)=3Q(T,y)\in 3\mathbb Z$. This means that $m_0\in A_0$ is a long root of $L_0$.\\

Let $det((\mathbb Z m_0 + \mathbb Z h^2)^{\wedge})=6$. By Claim 3, in this case $\mathbb Z m_0 + \mathbb Z h^2$ is saturated in $L$ in this case. Hence $Q(m_0,m_0)=2$ by formula (3), i.e. $m_0 \in A_0$ is a root.\\

This proves that (b) implies (a).\\ 

Now let's prove that (a) implies (b).\\

Let $m_0 \in A_0$ be a root, i.e. $Q(m_0,m_0)=2$. Then $\mathbb Z m_0 \subset L$ is saturated, and by (3) and (4) $\mathbb Z m_0 + \mathbb Z h^2 \subset L$ is saturated as well, and $det((\mathbb Z m_0 + \mathbb Z h^2)^{\wedge})=6$. So, statement (b) holds in this case.\\

Let $m_0 \in A_0$ be a long root of $L_0$. In particular, $Q(m_0,m_0)=6$.\\

By Claim 4, $\mathbb Z m_0 + \mathbb Z h^2$ is not saturated in $L$. Hence by formula (4), we get $det((\mathbb Z m_0 + \mathbb Z h^2)^{\wedge})=2$. So, statement (b) holds in this case too.\\

Consequently, statement (a) implies statement (b). {\it QED}\\

{\bf Lemma 3.4} {\it Let $W={\mathbb R}^{21-\rho} \oplus {\mathbb R}^2$ be an $\mathbb R$-vector space with quadratic form $Q(z,w)=z_1^2+...+z_{21-\rho}^2-w_1^2-w_2^2$. Let $A_1 \subset W$ be a discrete (hence closed and free) subgroup of rank $21-\rho +2$ (so that it forms a complete lattice). Let $P \subset W {\otimes}_{\mathbb R} {\mathbb C}$ be a union of a countable number of hyperplanes. \\
Then there exists an $\mathbb R$-vector subspace $\Sigma \subset W$, such that $sign(Q|_{\Sigma})=(21-\rho,0)$, $\Sigma \cap A_1=0$. Moreover, one can choose $\Sigma$ in such a way that there exists $f_0\in W\otimes_{\mathbb R}\mathbb C$ such that $f_0\neq 0$, $Q(f_0,\Sigma)=0$, $f_0\notin P$, $Q(f_0,f_0)=0$ and $Q(f_0,\bar f_0)<0$.}\\

{\it Proof:} Note that the condition $Q(f_0,\bar f_0)<0$ follows automatically from the others. Condition $f_0\neq 0$ is also vacuous.\\

Take $\Sigma ' =\mathbb R^{21-\rho}\subset {\mathbb R}^{21-\rho} \oplus {\mathbb R}^2=W$. Now take $\Sigma$ to be a small generic perturbation of $\Sigma '$. {\it QED}\\ 

In the next lemma and its proof we will use the notation of Corollary 3.2.\\

{\bf Lemma 3.5} {\it There exists an $\mathbb R$-vector subspace $\Sigma \subset L_0 \otimes_{\mathbb Z} \mathbb R$, such that $sign(Q|_{\Sigma})=(20,0)$, $\Sigma \cap L_0=A_0$, and there exists $f_0\in L_0 \otimes_{\mathbb Z}\mathbb C$ such that $f_0\neq 0$, $f_0 \in im(\tau)$ and $Q(f_0,\Sigma)=0$.}\\

{\it Proof:} Let $W=(A_0 \otimes_{\mathbb Z}\mathbb R)^{\bot} \subset L_0 \otimes_{\mathbb Z}\mathbb R$. Then $sign(Q|_W)=(21-\rho,2)$.\\

Let $A_1=W \cap L_0$ and $P\subset W \otimes_{\mathbb R}\mathbb C$ be the intersection of $W \otimes_{\mathbb R}\mathbb C$ with the union of hyperplanes $D_M$ from Theorem 3.1. This is well defined provided $A_0$ contains neither roots, nor long roots of $L_0$ by Lemma 3.3.\\

By Lemma 3.4, there exists an $\mathbb R$-vector subspace $\Sigma ' \subset W$, such that $sign(Q|_{\Sigma '})=(21-\rho,0)$, $\Sigma ' \cap A_1=0$, and there exists $f_0\in W \otimes_{\mathbb R}\mathbb C \subset L_0 \otimes_{\mathbb Z}\mathbb C$, such that $f_0 \neq 0$, $Q(f_0, \Sigma ')=0$, $f_0 \notin P$, $Q(f_0,f_0)=0$ and $Q(f0,\bar f_0)<0$.\\

Now it is enough to take $\Sigma = A_0 \otimes_{\mathbb Z}\mathbb R + \Sigma ' \subset L_0 \otimes_{\mathbb Z}\mathbb R$. {\it QED}\\ 

{\it Proof of Corollary 3.2:}\\

Choose $\Sigma \subset L_0 \otimes_{\mathbb Z} \mathbb R$ as in Lemma 3.5. Then there exists a cubic fourfold $X$ and an isomorphism of lattices $H^4(X,\mathbb Z)\cong L$, which identifies the square of the hyperplane section class with $h^2$, and $H^{3,1}(X)$ with $\mathbb C f_0 \subset L_0 \otimes_{\mathbb Z}\mathbb C$.\\

We then get the following identifications: 
\begin{itemize}
\item $H^{3,1}(X)=\mathbb C f_0$
\item $H^{1,3}(X)=\mathbb C \bar{f_0}$
\item $H^{2,2}(X)=\Sigma {\otimes}_{\mathbb R} {\mathbb C}$
\end{itemize}

In particular, $A_0(X)=H^4(X,\mathbb Z)_0\cap H^{2,2}(X)=(\Sigma {\otimes}_{\mathbb R} {\mathbb C})\cap L_0=\Sigma \cap L_0=A_0$. {\it QED}\\

{\bf Corollary 3.6} {\it Let $T\subset L_0$ be a primitive sublattice of signature $(21-\rho,2)$. Then the existence of a cubic fourfold $X$ and an isomorphism of lattices $H^4(X,\mathbb Z)\cong L$, identifying the self-intersection of the hyperplane section class with $h^2$, and $T(X)$ with $T$, is equivalent to the condition that the orthogonal complement $A_0=T_{L_0}^{\bot}$ of $T$ in $L_0$ contains neither roots, nor long roots of $L_0$.}\\ 

{\it Proof:} This follows from Corollary 3.2 by taking orthogonal complements. {\it QED}\\

\section{Existence of primitive embeddings.}

Let $2 \leq \rho \leq 21$. We are going to study even lattices of signature $(\rho-1,0)$, which admit a primitive embedding into $L_0$.\\

For the purpose of reference, let us state some simple properties of lattice $L_0$.\\ 

{\bf Lemma 4.1}  
{\it \begin{enumerate}
\item[(a)] $sign(L_0)=(20,2)$;
\item[(b)] $sign(q_{L_0})\equiv 2(mod 8)$;
\item[(c)] $A_{L_0}=A_{\begin{pmatrix} 2 & 1\\ 1 & 2 \end{pmatrix}}={\mathbb Z} / 3 {\mathbb Z}$;
\item[(d)] $q_{L_0} \colon {\mathbb Z} / 3 {\mathbb Z} \longrightarrow \mathbb Q / 2 \mathbb Z$, $1 \mapsto \frac{2}{3} + 2 \mathbb Z$;
\item[(e)] $l(A_{L_0})=1$.
\end{enumerate} }

We will also use the following observations:\\

{\bf Lemma 4.2} {\it Let $S\longrightarrow L$ be a primitive embedding of lattices, $M=L/S$. Then we have the following exact sequences of abelian groups:
\begin{itemize}
\item $L/S \longrightarrow S^{*}/S \longrightarrow S^{*}/L \longrightarrow 0$;
\item $0 \longrightarrow (L+M^{*})/L \longrightarrow L^{*}/L \longrightarrow S^{*}/L \longrightarrow 0$.
\end{itemize} }

{\bf Corollary 4.3} {\it $l(A_S)\leq rk(L/S)+l(S^{*}/L)\leq rk(L/S)+l(A_L)$.}\\

{\bf Remark} As we have already mentioned above, $l(A_S)\leq rk(S)$ for any lattice $S$.\\

{\bf Lemma 4.4} {\it Lattice $L_0$ is the unique even lattice with signature $(20,2)$ and discriminant-form $q_{L_0} \colon {\mathbb Z} / 3 {\mathbb Z} \longrightarrow \mathbb Q / 2 \mathbb Z$, $1 \mapsto \frac{2}{3} + 2 \mathbb Z$.}\\ 

{\it Proof:} This follows from Corollary 1.13.3 of \cite{Nikulin}. {\it QED}\\

{\bf Lemma 4.5} {\it The existence of primitive embedding of $A_0$ into $L_0$ depends only on the genus of $A_0$, i.e. on the signature $sign(A_0)=(\rho-1,0)$ and the discriminat-form $q_{A_0} \colon A_{A_0} \longrightarrow \mathbb Q / 2 \mathbb Z$.}\\ 

{\it Proof:} This follows from Lemma 2.3 of \cite{Morrison}. {\it QED}\\ 

{\bf Lemma 4.6} {\it Let $\Omega = U^{22}$ (the unique even unimodular lattice of signature $sign(\Omega)=(22,22)$), $T=L_0(-1)$, $sign(T)=(2,20)$, $l(A_T)=1$, $q_T=-q_{L_0}$ ($A_T=A_{L_0}$). Then:
\begin{itemize}
\item[(a)] $T=L_0(-1)$ is the unique even lattice in its genus;
\item[(b)] There exists a unique primitive embedding $L_0 \hookrightarrow \Omega$; 
\item[(c)] There exists a unique primitive embedding $L_0(-1) \hookrightarrow \Omega$;
\item[(d)] There exists a unique primitive embedding $A_0 \hookrightarrow \Omega$;
\item[(e)] $L_0 \bot L_0(-1)$ in $\Omega$.
\end{itemize}}

{\it Proof:} (a) follows from Corollary 1.13.3 of \cite{Nikulin}. (b), (c) and (d) follow from Corollary 1.12.3 and Corollary 1.14.4 of \cite{Nikulin}. (e) follows from Corollary 1.6.3 of \cite{Nikulin}. {\it QED}\\ 

In the next lemma we use the same idea as was used in Proposition 1.15.1 by \cite{Nikulin}.\\

{\bf Lemma 4.7} {\it Let $S$ be any (nonzero) even lattice. Assume that $sign(S) \leq sign(L_0)$. Consider the following statements:
\begin{itemize}
\item[(a)] There exists a primitive embedding $S \hookrightarrow L_0$.
\item[(b)] There exists a primitive embedding of even lattices $S \hookrightarrow V$, such that $S_V^{\bot} \cong L_0(-1)$, and $V$ admits a primitive embedding into $\Omega$.
\item[(c)] There exists a primitive embedding $S \oplus L_0(-1) \hookrightarrow \Omega = U^{22}$. 
\item[(d)] There exists $h\in A_S$ an element of order $3$, such that $q_S(h)=\frac{2}{3} \in \mathbb Q / 2 \mathbb Z$, such that if $V$ is (the unique) even lattice with signature $sign(V)=sign(S)+(2,20)$ and discriminant-form $q_V=q_S|_{A_V} \colon A_V \rightarrow \mathbb Q / 2 \mathbb Z$, where $A_V=\{  e, h, h^2 \}_{A_S}^{\bot} \subset A_S$, then there exists a primitive embedding $V \hookrightarrow \Omega$.
\end{itemize}
Then: (a) holds if and only if (b) holds, and (b) holds if and only if either (c) or (d) holds.}\\ 

{\bf Remark.} By Theorem 1.12.2 in \cite{Nikulin}, a primitive embedding $V \hookrightarrow \Omega$ in the situation of (d) exists if and only if there exists an even lattice with signature equal to $sign(L_0)-sign(S)$ and discriminant-form $-q_V$. This is satisfied, when $rk(L_0)-rk(S) \geq l(A_V)+1$.\\

{\bf Remark.} Note that $A_S=\{ e, h, h^2 \}\oplus A_V$, and so $l(A_V)\leq l(A_S)$.\\

{\it Proof:} Suppose we have a primitive embedding $S \hookrightarrow L_0$. Compose it with the primitive embedding $L_0 \hookrightarrow \Omega$ from Lemma 4.6 (b) to get a primitive embedding $S \hookrightarrow \Omega$, and hence an embedding of lattices $S \oplus L_0(-1) \hookrightarrow \Omega$. Take $V$ to be the saturation of the image of $S \oplus L_0(-1)$ in $\Omega$. Hence (a) implies (b).\\

Suppose (b) holds. Take the composition of primitive embeddings of lattices $S \hookrightarrow V$ and $V \hookrightarrow \Omega$ to get a primitive embedding $S \hookrightarrow \Omega$. Since the orthogonal complement of the image of $S$ in $\Omega$ contains $L_0(-1)$, the primitive embedding $S \hookrightarrow \Omega$ factors trough a primitive embedding $S \hookrightarrow L_0$.\\  

This establishes equivalence of (a) and (b).\\

By Proposition 1.5.1 \cite{Nikulin} the existence of a primitive embedding of $S$ into an even lattice $V$ with discriminant-form $q$, such that the orthogonal complement of $S$ in $V$ is isomorphic to $L_0(-1)$, is equivalent to the existence of a subgroup $H \subset A_{L_0(-1)}$ that admits an embedding $\gamma \colon H \hookrightarrow A_S$, such that $q_S \circ \gamma = -q_{L_0(-1)}| H$, and $(q_{L_0(-1)}\oplus q_S| (\Gamma_{\gamma})^{\bot})/\Gamma_{\gamma}\cong q$, where $\Gamma_{\gamma}$ is the graph of $\gamma$ in $A_{L_0(-1)}\oplus A_S$.\\

So, the part of (b) about the existence of the primitive embedding of $S$ can be restated as follows. (Recall that $q_{L_0(-1)}=-q_{L_0} \colon \mathbb Z / 3 \mathbb Z \rightarrow \mathbb Q / 2 \mathbb Z$, $1 \mapsto \frac{2}{3} + 2 \mathbb Z$.)\\

Either we can take $V=S \oplus L_0(-1)$ (which corresponds to $H=\{ e \}$), and in this case the condition is: 

\begin{center}
$V=S \oplus L_0(-1)$ admits a primitive embedding into $\Omega$,
\end{center}

or (if we take $H=A_{L_0(-1)}$):\\
\begin{center}
There exists an element $h \in A_S$ of order $3$, such that $q_S(h)=\frac{2}{3} + 2 \mathbb Z$. There exists a primitive embedding into $\Omega$ for some even lattice $V$ with signature $sign(V)=sign(S)+(2,20)$ and discriminant-form $q_V=q_S|_{A_V} \colon A_V \longrightarrow \mathbb Q / 2 \mathbb Z$, where $A_V=\{  e, h, h^2 \}_{A_S}^{\bot} \subset A_S$.
\end{center}

{\bf Remark.} By Corollary 1.13.3 in \cite{Nikulin} there exists a unique (up to an isomorphism) even lattice $V$ with signature and discriminant-form as above.\\

So, we proved that (b) holds if and only if either (c) or (d) holds. {\it QED}\\ 

{\bf Corollary 4.8} {\it Let $S$ be any even lattice. Assume that $sign(S) \leq sign(L_0)$. Then either of the following two conditions implies the existence of a primitive embedding $S \hookrightarrow L_0$:\\
\begin{enumerate}
\item $rk(L_0)-rk(S) \geq l(A_S)+2$. 
\item $rk(L_0)-rk(S) \geq l(A_S)+1$ and there exists $h\in A_S$ an element of order $3$, such that $q_S(h)=\frac{2}{3} \in \mathbb Q / 2 \mathbb Z$.
\end{enumerate}}
{\it Proof:} (1) follows from Lemma 4.7 and Corollary 1.12.3 in \cite{Nikulin}.\\

(2) follows from Lemma 4.7 and the remarks after it. {\it QED}\\ 

{\bf Theorem 4.9} {\it If $2\leq \rho \leq 11$, and $A_0$ is an even lattice of signature $sign(A_0)=(\rho - 1,0)$, then there exists a primitive embedding $A_0 \hookrightarrow L_0$.\\

Moreover, if $\rho \leq 10$, then the lattice $T=(A_0)_{L_0}^{\bot}$ is uniquely determined by its signature $(21-\rho, 2)$ and discriminant-form.}\\ 

{\it Proof:} We will use only the first condition from Corollary 4.8.\\
 
Recall that $sign(L_0)=(20,2)$, and so $rk(L_0)-rk(S)=23-\rho$. Since $l(A_S) \leq rk(S)=\rho - 1$, the first condition in Corollary 4.8 would follow from the inequality $23-\rho \geq \rho + 1$, i.e. $\rho \leq 11$.\\

This proves the first statement.\\

The second statement follows from \cite{Nikulin}, Corollary 1.13.3. Indeed, $rk(T)=rk(L_0)-rk(A_0)=23-\rho$, and $l(A_T)\leq rk(L_0/T)+l(A_{L_0})=rk(A_0)+1=\rho$ (by Corollary 4.3). Hence $rk(T)\geq l(A_T)+2$, when $\rho \leq 10$. This is the sufficient condition from \cite{Nikulin}, Corollary 1.13.3. {\it QED}\\ 

{\bf Theorem 4.10} {\it If $13\leq \rho \leq 21$, then for any even lattice $T$ of signature $sign(T)=(21-\rho,2)$, there exists a primitive embedding $T \hookrightarrow L_0$.}\\

{\it Proof:} We will use only the first condition from Corollary 4.8.\\

Since $rk(L_0)-rk(S)=\rho-1$, and $l(A_T)\leq rk(T)=23-\rho$, the first condition of Corollary 4.8 would follow from the inequality $\rho-1 \geq 25-\rho$, i.e. $\rho \geq 13$. {\it QED}\\ 

{\bf Lemma 4.11} {\it Let $S$ be an even lattice and consider a primitive embedding $S\hookrightarrow L_0$. 
\begin{enumerate}
\item[(a)] If the primitive embedding $S\hookrightarrow L_0$ is obtained from a primitive embedding $S\oplus L_0(-1)\hookrightarrow \Omega$ (case (c) of Lemma 4.7), then (the image of) $S$ does not contain long roots of $L_0$.
\item[(b)] If the primitive embedding $S\hookrightarrow L_0$ is obtained from a primitive embedding $S\oplus L_0(-1)\hookrightarrow \Omega$ (case (c) of Lemma 4.7), then the orthogonal complement $K$ of (the image of) $S$ contains a long root of $L_0$, if and only if there exists $\Delta \in K^{*}$ such that its class $\bar{\Delta}$ in $A_K$ has a form $(0,\pm 1)$ under the isomorphism $A_K\cong A_S\oplus A_{L_0}$ ($A_{L_0}\cong \mathbb Z/ 3 \mathbb Z$) corresponding to the primitive embedding $S\oplus L_0(-1)\hookrightarrow \Omega$ and $Q_{K^{*}}(\Delta,\Delta)=\frac{2}{3}$.
\item[(c)] If the primitive embedding $S\hookrightarrow L_0$ is obtained from a primitive embedding $V\hookrightarrow \Omega$ (case (d) of Lemma 4.7), then the orthogonal complement of (the image of) $S$ does not contain long roots of $L_0$.   
\item[(d)] If the primitive embedding $S\hookrightarrow L_0$ is obtained from a primitive embedding $V\hookrightarrow \Omega$ (case (d) of Lemma 4.7), then (the image of) $S$ contains a long root of $L_0$, if and only if there exists $\delta \in S$ such that $Q_S(\delta,\delta)=6$ and for any $\sigma \in S^{*}$ with $\bar{\sigma}\in A_V\subset A_S$ we have: ${\sigma}(\delta)\in 3\mathbb Z$. \\
\end{enumerate} }

{\it Proof:} (a) Let $K$ be an even lattice with invariants $sign(K)=sign(L_0)-sign(S)$, $A_K \cong A_S\oplus A_{L_0}$, $q_K\cong (-q_S)\oplus q_{L_0}$, which is the orthogonal complement of $S\oplus L_0(-1)$ in its primitive embedding into $\Omega$. Then $\Omega$ is the overlattice of $S\oplus L_0(-1)\oplus K$ corresponding to the subgroup $H\subset A_S \oplus A_{L_0(-1)}\oplus A_K$ that is the graph of the isomorphism $A_K \cong A_S\oplus A_{L_0}$ above (see Proposition 1.5.1 and Proposition 1.6.1 in \cite{Nikulin}). In other words, if we identify $A_K$ with $A_S\oplus A_{L_0}$ via this isomorphism (and $ A_{L_0}$ with $ A_{L_0(-1)}$), then 
$$
H=\{ (\bar s,\bar l,\bar k)\in A_S \oplus A_{L_0(-1)}\oplus A_K \quad | \quad \bar k=\bar s+\bar l \}
$$

The construction used in Proposition 1.5.1 in \cite{Nikulin} allows now to descibe explicitely this overlattice $\Omega$ in terms of the other lattices as follows: 
$$
\Omega=\{ (s,l,k)\in S^{*}\oplus L_0(-1)^{*}\oplus K^{*} \quad  |\quad  \bar k=\bar s+\bar l  \}
$$ 
where bar denotes the classes of elements in the corresponding quotients.\\ 

Note that $K=S_{L_0}^{\bot}$. \\

To simplify notation, let us denote $a=|A_S|$.\\

The intersection form on $\Omega$ is defined as follows:\\
$$
Q((s_1,l_1,k_1),(s_2,l_2,k_2))=\frac{1}{9 a^2} \cdot (Q_S(3a \cdot s_1,3a\cdot s_2)-Q_{L_0}(3a\cdot l_1,3a\cdot l_2)+Q_K(3a\cdot k_1,3a\cdot k_2))
$$
 
Since the orthogonal complement of $L_0(-1)$ in $\Omega$ (with respect to the given embedding $S\oplus L_0(-1)\oplus K \subset \Omega$) is isomorphic to $L_0$, for simplicity of notation, we denoted it by the same letter $L_0$.\\ 
 
As a sublattice of $\Omega$ this orthogonal complement can be described as follows (using the explicit description of the intersection form above): 
$$
L_0 = (L_0(-1))_{\Omega}^{\bot} = \{ (s,l,k)\in \Omega | l=0, \bar{k}=\bar{s}+0 \}
$$

Suppose there exists $\delta \in S$, such that $Q(\delta, \delta)=6$, and $Q(\delta, L_0) \subset 3 \mathbb Z$ (i.e., $\delta$ is a long root).\\

In terms of the bilinear form defined above this means that
$$
3\quad  \big | \quad \frac{1}{9a^2}\cdot Q(3a\cdot \delta,3a\cdot s)
$$
for each $s \in S^{*}$.\\

This implies that all elements of the set $Q_S(\delta, |A_S| \cdot S^{*})$ are integers divisible by $3 \cdot |A_S|$. In other words, for any $\sigma \in S^{*}$, $\sigma (\delta)$ is divisible by $3$. This is impossible.\\

Indeed, since $Q_S(\delta,\delta)=6$, it follows, that $\mathbb Z \delta \subset S$ is a primitive sublattice. Hence $\delta$ can be complemented to a basis in $S$: $\{ e_1=\delta, e_2, \ldots \}$. Now take $\sigma \in S^{*}$ to be the first element of the dual basis, i.e. $\sigma (e_1)=1$, and $\sigma (e_i)=0$ for any $i \geq 2$. In particular, we get that $\sigma (\delta)=1$, which is {\bf not} divisible by $3$. \\

So, we arrive at contradiction and conclude that there is no long root in this situation.\\

(b) Let $p \colon K^{*}\rightarrow A_{L_0}$ denote the projection, i.e. the composition of group homomorphisms $K^{*}\rightarrow A_K\cong A_S\oplus A_{L_0}\rightarrow A_{L_0}$.\\

Then $\delta \in K$ is a long root of $L_0$, if and only if $Q_K(\delta,\delta)=6$ and $(ker(p))(\delta)\subset 3\mathbb Z$. In particular, $Q_K(x,\delta)\in 3\mathbb Z$ for any $x\in K$.\\

Let $\delta \in K$ be a long root of $L_0$. Then we can take $\Delta=\frac{1}{3}\cdot \delta \in K^{*}$. If $\bar{\Delta}$ has a form $(\bar{\theta},\bar{\lambda})$ under the isomorphism $A_K\cong A_S\oplus A_{L_0}$, then the condition $(ker(p))(\delta)\subset 3\mathbb Z$ translates into saying that for any ${\bar{\theta}}_1\in A_S$,  $Q_{K^{*}}(({\bar{\theta}}_1,0),(\bar{\theta},\bar{\lambda}))\subset \mathbb Z$, i.e. $b_S({\bar{\theta}}_1,\bar{\theta})=0$ in $\mathbb Q / \mathbb Z$. Since $b_S$ is nondegenerate (and so $\bar{\theta}=0$) and $Q_{K^{*}}(\bar{\Delta},\bar{\Delta})=\frac{1}{9}\cdot Q_K(\delta,\delta)=\frac{2}{3}\notin 2\mathbb Z$, we conclude that $q_K(\bar{\lambda})=\frac{2}{3}+2\mathbb Z$, and so $\bar{\Delta}$ satisfies the required conditions.\\

Suppose that $\Delta\in K^{*}$ is such that its class $\bar{\Delta}\in A_K$ has a form $(0,\pm 1)$ and $Q_{K^{*}}(\Delta,\Delta)=\frac{2}{3}$. Then $\delta =3\Delta \in K$ will be a long root of $K$. \\
 
(c) According to Proposition 1.12.2 of \cite{Nikulin}, the existence of a primitive embedding $V\hookrightarrow \Omega$ amounts to the existence of an even lattice $K$ of signature $sign(K)=sign(L_0)-sign(S)$ and a group isomorphism $\gamma \colon A_K\cong A_V$ such that $q_S\circ \gamma=-q_K$. \\

Then $\Omega$ can be viewed as an overlattice of $V\oplus K$, which is the preimage of the graph ${\Gamma}_{\gamma}$ of $\gamma$ under the natural projection $V^{*}\oplus K^{*}\rightarrow A_V\oplus A_K$
$$
\Omega = \{ (v,\kappa)\in V^{*}\oplus K^{*} \mid \quad {\gamma}(\bar{\kappa})=\bar{v} \} 
$$

Note that according to the construction used in Proposition 1.5.1 of \cite{Nikulin}, we can describe $V$ explicitly as an abelian subgroup of $S^{*}\oplus L_0(-1)^{*}$ as well:
$$
V=\{ (\theta,\lambda)\in S^{*}\oplus L_0(-1)^{*} \mid \quad \bar{\theta} =\bar{\lambda} \cdot h \}
$$

Moreover, we have a sequence of embeddings of abelian groups $S\oplus L_0(-1)\hookrightarrow V\hookrightarrow V^{*}\hookrightarrow S^{*}\oplus L_0(-1)^{*}$, where $V^{*}$ is the preimage of the orthogonal complement of the graph ${\Gamma}_{\gamma}$ as above, and so it can be described explicitly as follows:
$$
V^{*}=\{ (\theta,\lambda)\in S^{*}\oplus L_0(-1)^{*} \mid \quad b_S(\bar{\theta} -\bar{\lambda} \cdot h,h)=0 \}
$$

Let us denote $G=(h)^{\bot}_{A_S}\subset A_S$.\\

Hence
$$
\Omega = \{ (\theta,\lambda,\kappa)\in S^{*}\oplus L_0(-1)^{*}\oplus K^{*} \mid \quad \bar{\theta}-\bar{\lambda}h\in G\subset A_S, \quad {\gamma}(\bar{\kappa})=\bar{\theta}-\bar{\lambda}h   \}
$$

Then $L_0$ as a sublattice of $\Omega$ (i.e. the orthogonal complement of $L_0(-1)$) can be described as follows:
$$
L_0 = (L_0(-1))_{\Omega}^{\bot} = \{ (s,0,k)\in \Omega |\quad \bar{s}\in G,\quad {\gamma}(\bar{k})=\bar{s}+0 \}
$$

If $\delta \in K$ is a long root of $L_0$, then as in part (a), we get a primitive nonzero sablattice $\mathbb Z \delta \subset K$ such that for any $\kappa \in K^{*}$, ${\kappa}(\delta)\in 3\mathbb Z$, which is impossible.\\

Hence $K$ does not contain long roots of $L_0$. \\

(d) This follows immediately from the explicit description of $L_0$ as a sublattice of $\Omega$ in (c).\\ 

{\it QED}\\

{\bf Corollary 4.12} {\it Let $2\leq \rho \leq 11$, $S$ be an even lattice of signature $sign(S)=(\rho - 1,0)$. Then there exists a primitive embedding $S \hookrightarrow L_0$, such that (the image of) $S$ does not contain long roots of $L_0$.}\\

{\it Proof:} In Theorem 4.9 we proved the existence of a primitive embedding $S \hookrightarrow L_0$ by checking that there exists a primitive embedding $S\oplus L_0(-1)\hookrightarrow \Omega=U^{22}$. Lemma 4.11 implies that (the image of) $S$ does not contain long roots of $L_0$. {\it QED}\\

{\bf Lemma 4.13} {\it Let $S$ be an even lattice with signature $sign(S)=(t_{+},t_{-})\leq sign(L_0)$ and discriminant-form $q\colon A\rightarrow \mathbb Q/2\mathbb Z$, $A=A_S=S^{*}/S$.\\ 

Then there exists a primitive embedding $S\oplus L_0(-1)\hookrightarrow \Omega$, if and only if one of the following conditions holds:  
\begin{enumerate}
\item[A1] $l(A)=l(A_3)\leq 20-rk(S)$ 
\item[A2] $l(A)=l(A_3)=21-rk(S)$ and $\left( \frac{(-1)^{rk(S)}\cdot ({\theta}_3)^{v_3}}{3} \right)=$
\item[A3] {$22-rk(S)=l(A)>l(A_3)$ and the following conditions hold: 
		\begin{itemize}
		\item If for a prime $p\neq 2,3$, $l(A)=l(A_p)$, then $\left( \frac{(-1)^{rk(S)}\cdot 3\cdot ({\theta}_p)^{v_p}}{p} \right)=1$.	
		\item If $l(A)=l(A_3)+1$, then $\left( \frac{(-1)^{rk(S)+1}\cdot ({\theta}_3)^{v_3}}{3} \right)=1$.	
		\item If $l(A)=l(A_2)$ and $q_2\neq q_{\theta}^{(2)}(2)\oplus q_2^{'}$, then $v_2=1$ (in this case also $rk(S)>l(A)$).\\
		\end{itemize} }
\end{enumerate} }

{\it Proof:} Let us denote $V=S\oplus L_0(-1)$. Then $rk(\Omega)-rk(V)=22-rk(S)$, $A_V\cong A\oplus \mathbb Z/3\mathbb Z$ and $q_V=q\oplus (-q_{L_0})$.\\

Hence for any $p\neq 3$, $(A_V)_p=A_p$ and $(q_V)_p=q_p$, while $(A_V)_3=A_3\oplus \mathbb Z/3\mathbb Z$ and $(q_V)_3=q_3\oplus q_{-1}^{(3)}(3)$, $discr(K(q_V)_3))=-3\cdot discr(K(q_3))$. \\

In particular, $\mid A_V \mid=3\cdot \mid A \mid$ and either $l(A_V)=l(A)$ (if $l(A_3)<l(A)$), or $l(A_V)=l(A)+1$ (if $l(A_3)=l(A)$).\\

From Corollary 1.12.3 of \cite{Nikulin} we know that if $22-rk(S)\geq l(A_V)+1$, then a primitive embedding $V \hookrightarrow \Omega$ exists. So, we may assume that $l(A_V)=22-rk(S)$.\\

In this case, by Theorem 1.12.2 of \cite{Nikulin} a primitive embedding $V\hookrightarrow \Omega$ exists, if and only if the following conditions hold:
\begin{enumerate}
\item[(1)] If for a prime $p\neq 2,3$, $l(A_p)=l(A_V)$, then $\frac{(-1)^{t_{+}}\cdot 3\cdot \mid A \mid}{discr(K(q_p))}\in ({\mathbb Z}^{*}_p)^2 $.
\item[(2)] If $l(A_3)+1=l(A_V)$, then $\frac{(-1)^{t_{+}+1}\cdot \mid A \mid}{discr(K(q_3))}\in ({\mathbb Z}^{*}_3)^2 $.
\item[(3)] If $l(A_2)=l(A_V)$ and $q_2\neq q_{\theta}^{(2)}(2)\oplus q_2^{'}$, then $\pm \frac{3\mid A \mid}{discr(K(q_2))}\in ({\mathbb Z}^{*}_2)^2 $.
\end{enumerate}

Using Remark 2.1, we can rewrite these conditions as follows:
\begin{enumerate}
\item[(1)] If for a prime $p\neq 2,3$, $l(A_p)=l(A_V)$, then $\left( \frac{(-1)^{rk(S)}\cdot 3\cdot ({\theta}_p)^{v_p}}{p} \right)=1$.
\item[(2)] If $l(A_3)+1=l(A_V)$, then $\left( \frac{(-1)^{rk(S)+1}\cdot ({\theta}_3)^{v_3}}{3} \right)=1$.
\item[(3)] If $l(A_2)=l(A_V)$ and $q_2\neq q_{\theta}^{(2)}(2)\oplus q_2^{'}$, then $v_2=1$.
\end{enumerate}

This gives the lemma. {\it QED}\\

{\bf Lemma 4.14} {\it Let $S$ be an even lattice with signature $sign(S)=(t_{+},t_{-})\leq sign(L_0)$ and discriminant-form $q\colon A\rightarrow \mathbb Q/2\mathbb Z$, $A=A_S=S^{*}/S$.\\ 

Let $\tau \in A$ be an element of order 3 such that $q_S(\tau)=\frac{2}{3}\in \mathbb Q / 2\mathbb Z$ and $V$ be the even lattice with signature $sign(V)=sign(S)+(2,20)$ and discriminant-form $q_V=q_S {\mid}_{A_V}\colon A_V \rightarrow \mathbb Q / 2 \mathbb Z$, where $A_V=(\tau)^{\bot}_{A}\subset A$.\\

Then there exists a primitive embedding $V\hookrightarrow \Omega$, if and only if one of the following conditions holds:  
\begin{enumerate}
\item[B1] $l(A)=l(A_3)\leq 22-rk(S)$ 
\item[B2] $l(A)=l(A_3)=23-rk(S)$ and the following conditions hold: 
		\begin{itemize}
		\item If for a prime $p\neq 2,3$, $l(A)=1+l(A_p)$, then $\left( \frac{(-1)^{rk(S)}\cdot 3\cdot ({\theta}_p)^{v_p}}{p} \right)=1$.	
		\item If $l(A)=l(A_3)$, then $\left( \frac{(-1)^{rk(S)}\cdot ({\theta}_3)^{v_3}}{3} \right)=1$.	
		\item If $l(A)=1+l(A_2)$ and $q_2\neq q_{\theta}^{(2)}(2)\oplus q_2^{'}$, then $v_2=1$ (in this case also $rk(S)>l(A)$).	
		\end{itemize} 
\item[B3] $22-rk(S)=l(A)>l(A_3)$ and the following conditions hold: 
		\begin{itemize}
		\item If for a prime $p\neq 2,3$, $l(A)=l(A_p)$, then $\left( \frac{(-1)^{rk(S)}\cdot 3\cdot ({\theta}_p)^{v_p}}{p} \right)=1$.	
		\item If $l(A)=l(A_2)$ and $q_2\neq q_{\theta}^{(2)}(2)\oplus q_2^{'}$, then $v_2=1$ (in this case also $rk(S)>l(A)$).\\
		\end{itemize} 
\end{enumerate}}

{\it Proof:} Note that $rk(\Omega)-rk(V)=22-rk(S)$.\\

$A=A_V\oplus \mathbb Z / 3\mathbb Z$, and so for any $p\neq 3$, $A_p=(A_V)_p$ and $(q_V)_p=q_p$, while $A_3=(A_V)_3\oplus \mathbb Z/3\mathbb Z$ and $q_3=(q_V)_3\oplus q_1^{(3)}(3)$, $discr(K(q_3))=3\cdot discr(K(q_V)_3)$. \\

In particular, $\mid A_V \mid=\frac{\mid A \mid}{3}$ and either $l(A_V)=l(A)$ (if $l(A_3)<l(A)$), or $l(A_V)=l(A)-1$ (if $l(A_3)=l(A)$).\\

From Corollary 1.12.3 of \cite{Nikulin} we know that if $22-rk(S)\geq l(A_V)+1$, then a primitive embedding $V \hookrightarrow \Omega$ exists. So, we may assume that $l(A_V)=22-rk(S)$.\\

In this case, by Theorem 1.12.2 of \cite{Nikulin} a primitive embedding $V\hookrightarrow \Omega$ exists, if and only if the following conditions hold:
\begin{enumerate}
\item[(1)] If for a prime $p\neq 2,3$, $l(A_p)=l(A_V)$, then $\frac{(-1)^{t_{+}}\cdot \mid A \mid}{3\cdot discr(K(q_p))}\in ({\mathbb Z}^{*}_p)^2 $.
\item[(2)] If $l(A_3)-1=l(A_V)$, then $\frac{(-1)^{t_{+}}\cdot \mid A \mid}{discr(K(q_3))}\in ({\mathbb Z}^{*}_3)^2 $.
\item[(3)] If $l(A_2)=l(A_V)$ and $q_2\neq q_{\theta}^{(2)}(2)\oplus q_2^{'}$, then $\pm \frac{\mid A \mid}{3\cdot discr(K(q_2))}\in ({\mathbb Z}^{*}_2)^2 $.
\end{enumerate}

Using Remark 2.1, we can rewrite these conditions as follows:
\begin{enumerate}
\item[(1)] If for a prime $p\neq 2,3$, $l(A_p)=l(A_V)$, then $\left( \frac{(-1)^{rk(S)}\cdot 3\cdot ({\theta}_p)^{v_p}}{p} \right)=1$.
\item[(2)] If $l(A_3)-1=l(A_V)$, then $\left( \frac{(-1)^{rk(S)}\cdot ({\theta}_3)^{v_3}}{3} \right)=1$.
\item[(3)] If $l(A_2)=l(A_V)$ and $q_2\neq q_{\theta}^{(2)}(2)\oplus q_2^{'}$, then $v_2=1$.
\end{enumerate}

This gives the lemma. {\it QED}\\

{\bf Lemma 4.15} {\it Let $T$ be an even lattice with signature $sign(T)\leq sign(L_0)$ and discriminant-form $q_T\colon A_T \rightarrow \mathbb Q/2\mathbb Z$. Then there exists a primitive embedding $T\hookrightarrow L_0$ such that the orthogonal complement of the image of $T$ contains neither roots, nor long roots of $L_0$, if and only if one of the following conditions holds:
\begin{enumerate}
\item[A] There exists $\tau \in A_T$ an element of order $3$ such that $q_T(\tau)=\frac{2}{3}+2\mathbb Z\in \mathbb Q/2\mathbb Z$, and if we denote $G=(\tau)^{\bot}_{A_T}\subset A_T$, then there exists an even lattice $K$ of signature $sign(K)=sign(L_0)-sign(T)$ such that $A_K=G\subset A_T$, $q_K=-q_T{\mid}_G $, and $K$ has no roots. 
\item[B] There exists an even lattice $K$ without roots with signature $sign(K)=sing(L_0)-sign(T)$, $A_K=A_T\oplus A_{L_0}$ and discriminant-form $q_K=(-q_T)\oplus q_{L_0}$ such that there exists a group isomorphism $A_K\cong A_K$ (preserving $q_K$) such that if $p\colon K^{*}\rightarrow A_K\cong A_K$ is the composition, then there is no $\Delta \in p^{-1}((0,\pm 1))$ with $Q_{K^{*}}(\Delta,\Delta)=\frac{2}{3}$. 
\end{enumerate}}

{\it Proof:} This lemma follows from Lemma 4.7 and Lemma 4.11. Condition (A) describes situation in case (d) of Lemma 4.7. In this case $K=(T)^{\bot}_{L_0}$ never contains long roots of $L_0$ by Lemma 4.11 (c).\\

Condition (B) corresponds to the case (c) of Lemma 4.7. A criterion of the existence of long roots in this situation is given in Lemma 4.11 (b). {\it QED}\\

\section{Conclusions.}

Now we combine together some of the results of the previous two sections in the following Corollaries.\\

{\bf Corollary 5.1} {\it Let $2\leq \rho \leq 11$, and $A_0$ be an even lattice of signature $(\rho-1,0)$, which does not contain roots. Then there exists a cubic fourfold $X$ and an isomorphism of lattices $A_0(X) \cong A_0$.}\\

{\it Proof:} By Corollary 4.12 there exists a primitive embedding $A_0 \hookrightarrow L_0$, such that the image of $A_0$ does not contain long roots of $L_0$.\\ 

Hence we can view $A_0$ as a primitive sublattice of $L_0$ with no roots and no long roots of $L_0$. So, by Corollary 3.2 there exists a cubic fourfold $X$ with the required property. {\it QED}\\

{\bf Corollary 5.2} {\it Let $13\leq \rho \leq 21$, and $T$ be an even lattice of signature $(21 - \rho,2)$. Suppose one can choose a primitive embedding of $T$ into $L_0$ (existence of such embeddings was checked in Theorem 4.10) in such a way that its orthogonal complement contains neither roots, nor long roots of $L_0$. Then there is a cubic fourfold $X$ and an isomorphism of lattices $T(X)  \cong T$.}\\ 

{\it Proof:} This is a direct consequence of Corollary 3.6. {\it QED}\\

{\bf Theorem 5.3} {\it Let $S$ be an even lattice without roots, $sign(S)\leq sign(L_0)$. Then there exists a primitive embedding $S \hookrightarrow L_0$ such that the image of $S$ does not contain long roots of $L_0$, if and only if one of the following conditions holds:
\begin{enumerate}
				\item[A1] $l(A)=l(A_3)\leq 20-rk(S)$ 
				\item[A2] $l(A)=l(A_3)=21-rk(S)$ and $\left( \frac{(-1)^{rk(S)}\cdot ({\theta}_3)^{v_3}}{3} \right)=1$
				\item[A3] {$22-rk(S)=l(A)>l(A_3)$ and the following conditions hold: 
						\begin{itemize}
						\item If for a prime $p\neq 2,3$, $l(A)=l(A_p)$, then $\left( \frac{(-1)^{rk(S)}\cdot 3\cdot ({\theta}_p)^{v_p}}{p} \right)=1$.	
						\item If $l(A)=l(A_3)+1$, then $\left( \frac{(-1)^{rk(S)+1}\cdot ({\theta}_3)^{v_3}}{3} \right)=1$.	
						\item If $l(A)=l(A_2)$ and $q_2\neq q_{\theta}^{(2)}(2)\oplus q_2^{'}$, then $v_2=1$ (in this case also $rk(S)>l(A)$).\\
						\end{itemize} }
\item[B] There exists $h\in A$ an element of order $3$ such that $q_S(h)=\frac{2}{3}\in \mathbb Q/2 \mathbb Z$, there is no $\delta \in S$ with the property that $Q_S(\delta, \delta)=6$ and $\sigma (\delta)\in 3 \mathbb Z$ for any $\sigma \in S^{*}$ such that $Q_{S^{*}}(\sigma, h)\in \mathbb Z$, and one of the following conditions holds: 
				\begin{enumerate}
				\item[B1] $l(A)=l(A_3)\leq 22-rk(S)$ 
				\item[B2] $l(A)=l(A_3)=23-rk(S)$ and the following conditions hold: 
						\begin{itemize}
						\item If for a prime $p\neq 2,3$, $l(A)=1+l(A_p)$, then $\left( \frac{(-1)^{rk(S)}\cdot 3\cdot ({\theta}_p)^{v_p}}{p} \right)=1$.	
						\item If $l(A)=l(A_3)$, then $\left( \frac{(-1)^{rk(S)}\cdot ({\theta}_3)^{v_3}}{3} \right)=1$.	
						\item If $l(A)=1+l(A_2)$ and $q_2\neq q_{\theta}^{(2)}(2)\oplus q_2^{'}$, then $v_2=1$ (in this case also $rk(S)>l(A)$).	
						\end{itemize} 
				\item[B3] $22-rk(S)=l(A)>l(A_3)$ and the following conditions hold: 
						\begin{itemize}
						\item If for a prime $p\neq 2,3$, $l(A)=l(A_p)$, then $\left( \frac{(-1)^{rk(S)}\cdot 3\cdot ({\theta}_p)^{v_p}}{p} \right)=1$.	
						\item If $l(A)=l(A_2)$ and $q_2\neq q_{\theta}^{(2)}(2)\oplus q_2^{'}$, then $v_2=1$ (in this case also $rk(S)>l(A)$).\\
						\end{itemize} 
				\end{enumerate}
\end{enumerate}}

{\it Proof:} We combine together Lemma 4.7, Lemma 4.11, Lemma 4.13 and Lemma 4.14. {\it QED}\\

{\bf Remark.} Combined with Corollary 3.2, this theorem gives a classification of all possible lattices, which can appear as primitive saturated algebraic sublattices of intersection lattices of cubic fourfolds $A_0(X)$ (Corollary 5.3 in Section 1.2).\\

{\bf Theorem 5.4} {\it Let $2\leq \rho \leq 21$, $T$ be an even lattice with signature $sign(T)=(21-\rho,2)$ and discriminant-form $q_T\colon A_T \rightarrow \mathbb Q/ 2\mathbb Z$. Then there exists a cubic fourfold $X$ and an isomorphism of lattices $T(X)\cong T$, if and only if one of the following conditions holds:
\begin{enumerate}
\item[A] There exists $\tau\in A_T$ an element of order $3$ such that $q_T(\tau)=\frac{2}{3}+2\mathbb Z\in \mathbb Q / 2 \mathbb Z$, and if we denote $G=(\tau)^{\bot}_{A_T}\subset A_T$, then there exists an even lattice $K$ of signature $(\rho-1,0)$ such that $A_K=G\subset A_T$, $q_K=-q_T {\mid}_{G}$, and $K$ has no roots.
\item[B] There exists a primitive embedding of lattices $U^{23-\rho} \hookrightarrow T(-1)\oplus L_0$ such that the orthogonal complement $K$ of its image has no roots, and for any $\Delta \in (T(-1)\oplus {L_0}^{*})  \backslash (T(-1)\oplus {L_0})$ such that ${\Delta}(U^{23-\rho})=0$, ${\Delta}(3{\Delta})=\frac{1}{3}Q_{T(-1)\oplus {L_0}}(3{\Delta},3{\Delta})\neq 2$.
\end{enumerate}}

{\it Proof:} We apply Corollary 3.6 and Lemma 4.15.\\

Since isomorphisms $A_K=A_{T(-1)}\oplus A_{L_0} \cong A_{T(-1)}\oplus A_{L_0}$ (preserving discriminant-forms) correspond to the choices of an element $\xi \in A_T\oplus A_{L_0}$ of order $3$ with $q_K(\xi)=\frac{2}{3}+2\mathbb Z$ and a group isomorphism $(\xi)^{\bot}_{A_K}\cong A_{T(-1)}$ (preserving discriminant-forms), we can rephrase condition (B) in Lemma 4.15 as follows:
\begin{enumerate}
\item[B] There exists an even lattice $K$ without roots with signature $(\rho-1,0)$ such that there exists $\Delta \in K^{*}$ with the following properties:
\begin{itemize}
\item its class $\bar{\Delta}$ in $A_K=K^{*}/K$ has order $3$, $q_K(\bar{\Delta})=\frac{2}{3}+2\mathbb Z$, $(\bar{\Delta})^{\bot}_{A_K}\cong A_T$ and $q_K{\mid}_{(\bar{\Delta})^{\bot}_{A_K}}=-q_T$; 
\item ${\Delta}(3\Delta)+6\cdot \left(\pm {\Delta}(x)+\frac{Q_K(x,x)}{2} \right) \neq 2$ for any $x\in K$.
\end{itemize}
\end{enumerate}

Assuming $K$ exists, even lattices $U^{23-\rho}\oplus K$ and $T(-1)\oplus L_0$ have the same signature $(t_{+},t_{-})=(22,23-\rho)$ and discriminant-forms. We see that $t_{+}\geq 1$, $t_{-}\geq 1$ and $l(A_{T(-1)\oplus L_0})=l(A_T\oplus {\mathbb Z / 3\mathbb Z})\leq $  $ 1+l(A_T)\leq 1+rk(T)\leq $  $ 24-\rho \leq 45-\rho -2=t_{+}+t_{-}-2$. Hence these two lattices are isomorphic by Corollary 1.13.3 of \cite{Nikulin}.\\   

In addition, by Theorem 1.14.2 of \cite{Nikulin} every automorphism of $A_{T(-1)\oplus L_0}$ (preserving discriminant-forms) is induced by an automorphism of the lattice $T(-1)\oplus L_0$.\\ 

Hence condition (B) is equivalent to the condition stated in the theorem (the reader may want to recall Proposition 1.5.1 from \cite{Nikulin} at this point). {\it QED}\\ 

\section{Classification theorem for algebraic saturated sublattices of intersection lattices of cubic fourfolds.}

For a cubic fourfold $X$ we denote by $A(X)$ the saturated sublattice generated by its algebraic cycles:
$$
A(X)=H^4(X,\mathbb Z)\cap H^{2,2}(X)
$$

Note that $h^2 \in A(X)$ and $A_0(X)=(\mathbb Z h^2)_{A(X)}^{\bot}\subset A(X)$.\\

{\bf Theorem 6.1} {\it Let $A$ be an odd lattice of signature $sign(A)=(\rho,0)$, where $2\leq \rho \leq 21$. Then $A\cong A(X)$ for some cubic fourfold $X$ if and only if there exists $a\in A$, such that if we denote $A_0=(\mathbb Z a)_{A}^{\bot}$, then the following conditions are satisfied:
\begin{enumerate}
\item $Q_A(a,a)=3$.
\item $A_0$ is an even lattice.
\item There is {\it no} ${\delta}_0\in A_0$ such that $Q_A({\delta}_0,{\delta}_0)=2$.
\item There is {\it no} ${\delta}\in A_0$ such that $Q_A({\delta},{\delta})=6$ and for any $\alpha \in A^{*}$ with ${\alpha}(a)=0$ we have ${\alpha}(\delta)\in 3 \mathbb Z$.
\item For any $b\in A$ the integer $Q_A(a,b)^2-Q_A(b,b)$ is even.
\item There exists an even lattice $K$ of signature $sign(K)=sign(L)-sign(A)=(21-\rho,2)$, such that $A_A=A_K$ and $q_K \colon A_A=A^{*}/A \rightarrow \mathbb Q/ 2\mathbb Z$, $\alpha \mapsto \left[ ({\alpha}(a))^2-Q_{A^{*}}(\alpha,\alpha) \right] +2\mathbb Z$.
\end{enumerate}
}

{\bf Remark 6.2} Note that the bilinear form $b_K \colon A_K \otimes A_K \rightarrow \mathbb Q/ \mathbb Z$ of the quadratic form $q_K$ coincides with $-b_A$. In particular, it is nondegenerate. Together with condition $5$, this guarantees that the discriminant form $q_K$ given above is well-defined.\\ 

{\bf Remark 6.3} According to Theorem 1.10.1 from \cite{Nikulin}, condition $6$ holds if and only if the following conditions are satisfied (here for any prime $p$, $K((q_K)_p)$ denotes the unique $p-$adic lattice of rank $l((A^{*}/A)_p)$ with discriminant-form $(q_K)_p$, which exists according to Theorem 1.9.1 of \cite{Nikulin}):
\begin{enumerate}
\item[6.1] $19-\rho \equiv sign(q_K)\; mod(8)$.
\item[6.2] $23-\rho\geq l(A^{*}/A)$.
\item[6.3] If $p$ is an odd prime and $23-\rho=l((A^{*}/A)_p)$, then $\frac{discr(K((q_K)_p))}{|A^{*}/A|}\in ({\mathbb Z}_p^{*})^2$.
\item[6.4] If $23-\rho=l((A^{*}/A)_2)$ and $(q_K)_2\neq q_{\theta}^{(2)}(2)\oplus q_2^{'}$, then $\pm\frac{discr(K((q_K)_2))}{|A^{*}/A|}\in ({\mathbb Z}_2^{*})^2$.
\end{enumerate}

So, for condition $6$ to be satisfied, it is sufficient to require that $\rho \leq 11$ and $11-\rho \equiv sign(q_K) \; mod(8)$.\\

{\it Proof:} 

The full algebraic sublattice $A(X)$ of the intersection lattice $H^4(X,\mathbb Z)$ of a cubic fourfold $X$ is equal to the saturation of $\mathbb Z h^2 + A_0(X)$ in $H^4(X,\mathbb Z)$. So, we have to find a primitive embedding $A\hookrightarrow L$, such that some element $a \in A$ goes to $h^2\in L$, and its orthogonal complement becomes identified with $L_0$. As we remarked in section 2.1, it is sufficient to require that the image of $a\in A$ in $L$ has self-intersection $3$, and its orthogonal complement (in $L$) is even. (It follows then from the Nikulin's theory \cite{Nikulin} that there is an automorphism of $L$, which identifies $a$ with $h^2$, and its orthogonal complement with $L_0$.)\\

According to Corollary 3.2, the image of the orthogonal complement of $a$ in $A$ under this primitive embedding in addition should contain neither roots, not generalized roots of $L_0$.\\

So, we are looking for odd lattices $A$ that satisfy the following condition:
\begin{center}
There exists an even lattice $K$ and a primitive embedding $\phi \colon A\hookrightarrow L$, such that $Q_A(a,a)=3$, $(\mathbb Z \phi (a))_L^{\bot}$ is even and ${\phi}((\mathbb Z a)_{A}^{\bot})$ contains neither roots, nor long roots of $(\mathbb Z \phi (a))_L^{\bot}$.
\end{center}

From the version of Proposition 1.6.1 (\cite{Nikulin}) for odd lattices, we know that existence of primitive embeddings is related to the existence of lattices with specific properties. Using this criterion, the lattices $A$ we are looking for can be described as follows:
\begin{center}
$A$ is an odd lattice, such that there exists an even lattice $K$ with signature $sign(K)=sign(L)-sign(A)$ and a group isomorphism $\gamma \colon A_K \xrightarrow{\sim}  A^{*}/A=A_A$, such that $b_A\circ \gamma = -b_K$, and\\
there exists $a \in A$ such that $Q_A(a,a)=3$, $(\mathbb Z (a,0))_{{\pi}^{-1}({\Gamma}_{\gamma})}^{\bot}$ is even and $(\mathbb Z a)_{A}^{\bot}\oplus (0)\subset A\oplus K\subset {{\pi}^{-1}({\Gamma}_{\gamma})}$ contains neither roots, not long roots of $(\mathbb Z (a,0))_{{\pi}^{-1}({\Gamma}_{\gamma})}^{\bot}$.
\end{center}

Here we denote by ${\Gamma}_{\gamma}\subset A_A\oplus A_K$ the graph of the isomorphism $\gamma \colon A_K \xrightarrow{\sim} A_A$, and by $\pi \colon A^{*}\oplus K^{*}\rightarrow A_A\oplus A_K$ the projection map. Then by the construction, used in the proof of Proposition 1.6.1 of \cite{Nikulin}, we know how to find the overlattice $A\oplus K\subset L$ explicitly: it is isomorphic to the lattice ${{\pi}^{-1}({\Gamma}_{\gamma})}$: 
$$
A\oplus K \hookrightarrow {{\pi}^{-1}({\Gamma}_{\gamma})} \hookrightarrow A^{*}\oplus K^{*}.
$$

Since $|A_A|=|A_K|$, the intersection product on ${{\pi}^{-1}({\Gamma}_{\gamma})}$ is defined as follows:
$$
Q_{\infty}(({\alpha}_1,{\kappa}_1),({\alpha}_2,{\kappa}_2))=\frac{1}{|A_A|^2}\cdot \left[ Q_A(|A_A|\cdot {\alpha}_1,|A_A| {\alpha}_2)+Q_K(|A_K|\cdot {\kappa}_1,|A_K|\cdot {\kappa}_2) \right]=
$$

$$
=\frac{1}{|A_A|} \cdot \left[ {{\alpha}_2}(|A_A|\cdot {\alpha}_1) + {{\kappa}_2}(|A_K|\cdot {\kappa}_1) \right].
$$

From this we find, that lattice $(\mathbb Z (a,0))_{{\pi}^{-1}({\Gamma}_{\gamma})}^{\bot}$ inside ${{\pi}^{-1}({\Gamma}_{\gamma})}$ can be described explicitly as follows:
$$
(\mathbb Z (a,0))_{{\pi}^{-1}({\Gamma}_{\gamma})}^{\bot} = \{ (\alpha,\kappa)\in A^{*}\oplus K^{*} \quad | \quad {\gamma}(\bar{\kappa}) = \bar{\alpha} , {\alpha}(a)=0  \}.
$$

The condition that $(\mathbb Z a)_{A}^{\bot}\oplus (0)\subset A\oplus K \subset {{\pi}^{-1}({\Gamma}_{\gamma})}$ contains neither roots, nor long roots can be reformulated as follows:
\begin{center}
There does not exist ${\delta}_0\in (\mathbb Z a)_{A}^{\bot}$ such that $Q_A({\delta}_0,{\delta}_0)=2$, and\\
there does not exist $\delta \in (\mathbb Z a)_{A}^{\bot}$ such that $Q_A({\delta},{\delta})=6$ and for any $\xi \in (\mathbb Z (a,0))_{{\pi}^{-1}({\Gamma}_{\gamma})}^{\bot}$, $Q_{\infty}(\xi,\delta)\in 3\mathbb Z$\\
(i.e. for any $(\alpha,\kappa)\in A^{*}\oplus K^{*}$ such that ${\gamma}(\bar{\kappa})=\bar{\alpha}$ and ${\alpha}(a)=0$ we have ${\alpha}(\delta)\in 3 \mathbb Z$).
\end{center}

The condition that $(\mathbb Z (a,0))_{{\pi}^{-1}({\Gamma}_{\gamma})}^{\bot}$ is even can be restated as follows:
\begin{center}
For any $(\alpha,\kappa)\in A^{*}\oplus K^{*}$ such that ${\gamma}(\bar{\kappa})=\bar{\alpha}$ and ${\alpha}(a)=0$ we have $Q_{A^{*}}(\alpha,\alpha)+Q_{K^{*}}(\kappa,\kappa)  \in 2 \mathbb Z$.
\end{center}

We conclude that for an odd lattice $A$ of signature $sign(A)=(\rho,0)$ there exists a cubic fourfold $X$ and an isomorphism of lattices $A\cong A(X)$ if and only if there exists $a\in A$, such that if we denote $A_0=(\mathbb Z a)_{A}^{\bot}$, then the following conditions hold:
\begin{enumerate} 
\item $Q_A(a,a)=3$.
\item $A_0$ is an even lattice.
\item There is no ${\delta}_0 \in A_0$ such that $Q_A({\delta}_0,{\delta}_0)=2$.
\item There is no $\delta \in A_0$ such that $Q_A({\delta},{\delta})=6$ and for any $\alpha \in A^{*}$ with $\alpha (a)$ we have $\alpha(\delta) \in 3\mathbb Z$.
\item There exists an even lattice $K$ with signature $sign(K)=sign(L)-sign(A)$, such that there is a group isomorphism $\gamma \colon A_A \xrightarrow{\sim} A_K$ such that $b_K\circ \gamma =-b_A$ and for any $\alpha \in A^{*}$ with ${\alpha}(a)=0$ we have $Q_{A^{*}}(\alpha,\alpha)\equiv -q_K(\gamma (\bar{\alpha}))\; mod(2\mathbb Z)$.
\end{enumerate}

Taking into account that $A$ is odd, the last condition gives the description of the discriminant-form of the lattice $K$:
$$
q_K \colon A_A=A^{*}/A \rightarrow \mathbb Q /2\mathbb Z, \quad \quad \alpha \mapsto \left[ ({\alpha}(a))^2 - Q_{A^{*}}(\alpha,\alpha) \right] + 2\mathbb Z.
$$

In order for this discriminant-form to be well-defined we need to require that:
\begin{center}
for any $b\in A$, $Q_A(a,b)^2-Q_A(b,b)\in 2 \mathbb Z$.
\end{center}

Taking this all together, we obtain the criterion, stated in the theorem. {\it QED}\\

{\bf Remark.} If we consider the lattice with underlying free abelian group equal to $A$, but with the bilinear form given by $Q(x,y)=Q_A(x,a)\cdot Q_A(y,a)-Q_A(x,y)$, we will get a well-defined even lattice $A^{ev}$, such that $((A^{ev})^{*}/(A^{ev}))_p=(A^{*}/A)_p$ for any odd prime $p$, but $((A^{ev})^{*}/(A^{ev}))_2$ will be (at least, in general) different from $(A^{*}/A)_2$. In the latter case, $(A^{ev})^{*}/(A^{ev})$ will be a $\mathbb Z / 2\mathbb Z$-extension of $A^{*}/A$.\\

\section{Appendix. Existence of even lattices of given genus without roots 'locally'.}

From Theorem 1.9.1 of \cite{Nikulin} we know that for any nondegenerate finite quadratic form $q_p \colon A_{q_p} \rightarrow \mathbb Q_p /2 \mathbb Z_p$ there exists a unique p-adic lattice $K(q_p)$ of rank $l(A_{q_p})$, whose discriminant-form is isomorphic to $q_p$ (excluding the case, when $p=2$ and $q_2=q^{(2)}_{\theta}(2)\oplus q_2^{'}$). If $p=2$ and $q_2=q^{(2)}_{\theta}(2)\oplus q_2^{'}$ then there exist two such 2-adic lattices $K_{\alpha}(q_2)$, corresponding to two possible values of $\alpha$.\\ 

Observe that for any prime $p$, $\frac{discr(K(q_p))}{|A_{q_p}|}\in {\mathbb Z}^{*}_p$.\\

Suppose we are given non-negative integers $t_{+}$, $t_{-}$, $N=t_{+}+t_{-}\geq 1$, a finite abelian group $A=\oplus_{p}A_p$, where for each prime $p$, $A_p$ denotes a $p$-group, and for each prime $p$ a finite nondegenerate quadratic form $q_p \colon A_p \rightarrow \mathbb Q_p / 2 \mathbb Z_p$.\\

Let us introduce the following notation:
\begin{itemize}
\item for any odd prime $p$, ${\theta}_p=\frac{(-1)^{t_{-}}\cdot |A|}{discr(K(q_p))}\in {\mathbb Z}^{*}_{p}$;
\item if $q_2=q^{(2)}_{\theta}(2)\oplus q_2^{'}$, then let $K(q_2)$ denote the lattice $K_{\alpha}(q_2)$ such that $\frac{(-1)^{t_{-}+(N-l(A_2))/2}\cdot |A|}{discr(K_{\alpha}(q_2))}\in ({\mathbb Z}^{*}_{2})^2$, if such $\alpha$ exists, and any of the two lattices $K_{\alpha}(q_2)$ otherwise; 
\item ${\theta}_2 = \frac{(-1)^{t_{-}+(N-l(A_2))/2}\cdot |A|}{discr(K(q_2))}\in {\mathbb Z}^{*}_{2}$;
\item let $v_2=1$, if $5{\theta}_2\in ({\mathbb Z}^{*}_{2})^2$, and $0$ otherwise;
\item let ${\epsilon}_2=(-1)^{(N-l(A_2))(N-l(A_2)-2)/8+(N-l(A_2))(discr(K(q_2))/{\mid A_2 \mid}-1)/4}$, if $q_2=q^{(2)}_{\theta}(2)\oplus q_2^{'}$ or $N\geq l(A_2)+2$, and 1 otherwise.
\end{itemize}

{\bf Theorem 7.1} {\it With data and notation as above, there exists an even lattice $S$ with signature $(t_{+},t_{-})$ such that $q_S\cong \oplus_p q_p$, and for some prime $p$ the corresponding p-adic lattice $S_p$ has no roots, if and only if the following conditions are satisfied:
\begin{enumerate}
\item $N\geq l(A)$, $N\equiv l(A_2) \; (mod 2)$;
\item For any odd prime $p$ such that $N=l(A_p)$ we have that $\left( \frac{{\theta}_p}{p} \right)=1$;
\item If ${\theta}_2\notin ({\mathbb Z}^{*}_{2})^2$, then $q_2\neq q^{(2)}_{\theta}(2)\oplus q_2^{'}$, $N \geq 2+l(A_2)$ and $5{\theta}_2\in ({\mathbb Z}^{*}_{2})^2$;
\item ${\epsilon}_2=\prod\limits_{p} {\epsilon}_p(K(q_p))\cdot \prod\limits_{p\neq 2}{\left( \frac{{\theta}_p}{p} \right)}^{{log}_p(|A_p|)}\cdot (-1)^{t_{-}\cdot (t_{-}-1)/2+v_2\cdot(1+{log}_2|A_2|)}$ (here ${\epsilon}_p$ denotes Hasse invariant);
\item One of the following conditions holds:
			\begin{itemize}
			\item $v_2=0$, $N=l(A_2)$ and $K(q_2)$ has no roots;
			\item $\exists p\neq 2$ such that $N=l(A_p)$ and $K(q_p)$ contains no roots; 
			\item $\exists p\neq 2$ such that $N=l(A_p)+1$ and $K_{{\theta}_p}^{(p)}(1)\oplus K(q_p)$ contains no roots; 
			\item $\exists p\neq 2$ such that $N=l(A_p)+2$ and $K_{1}^{(p)}(1)\oplus K_{{\theta}_p}^{(p)}(1)\oplus K(q_p)$ contains no roots.  
			\end{itemize}
\end{enumerate}
}

{\it Proof:} According to Corollary 1.9.3 of \cite{Nikulin}, p-adic lattices corresponding to the lattice $S$ should have the form:
$$
S_p\cong K_{1}^{(p)}(1)^{t_p-1}\oplus K_{{\theta}_p}^{(p)}(1)\oplus K(q_p),
$$
where ${\theta}_p\in {\mathbb Z}_p^{*}/({\mathbb Z}_p^{*})^2$, $t_p=rk(S)-l(A_p)$, if $p$ is odd and $l(A_p)<rk(S)$,
$$
S_p\cong K(q_p),
$$
if $l(A_p)=rk(S)$, and 
$$
S_2\cong U^{(2)}(1)^{t_2-v_2}\oplus V^{(2)}(1)^{v_2}\oplus K_{\alpha}(q_2),
$$
where $0\leq v_2\leq 1$ and $t_2\geq v_2$.\\

We have $(-1)^{t_{-}} \mid A\mid =discr(S)=discr(S_2)$ $=(-1)^{t_2}\cdot (-3)^{v_2}\cdot discr(K_{\alpha}(q_2))$ $=discr(S_p)={\theta}_p\cdot discr(K(q_p))$, and $N=rk(S)=rk(S_2)=2t_2+l(A_2)$ $=rk(S_p)=t_p+l(A_p)$ for odd $p$. Then ${\theta}_p=(-1)^{t_{-}} \mid A \mid / discr(K(q_p))$ for odd $p$ (unless $l(A_p)=N$) and $(-1)^{t_{-}+(N-l(A_2))/2}\cdot \mid A \mid / discr(K_{\alpha}(q_2))=(-3)^{v_2}=5^{v_2}\in {\mathbb Z}_p^{*}/({\mathbb Z}_p^{*})^2$ (unless $l(A_2)=N$).\\

According to Theorem 1.10.3 of \cite{Nikulin}, p-adic lattices $S_p$ glue together and give a lattice $S$ over $\mathbb Z$, if and only if the following conditions hold:
\begin{enumerate}
\item[1.] For any odd prime $p$ such that $l(A_p)=N$, $(-1)^{t_{-}} \mid A \mid / discr(K(q_p))\in ({\mathbb Z}_p^{*})^2$.
\item[2.] If $N=l(A_2)$, then $(-1)^{t_{-}} \mid A \mid / discr(K_{\alpha}(q_2))\in ({\mathbb Z}_2^{*})^2$.
\item[2'.] If $N>l(A_2)$, then $(-3)^{v_2}\cdot (-1)^{t_{-}+(N-l(A_2))/2}\cdot \mid A \mid / discr(K_{\alpha}(q_2))\in ({\mathbb Z}_2^{*})^2$.
\item[3.] $(-1)^{t_{-}(t_{-}-1)/2}=\prod\limits_{p} {\epsilon}_p(S_p)$.
\end{enumerate}

Recalling standard formulas for the Hasse invariant (stated, for example, in section 10 of \cite{Nikulin}), we get ${\epsilon}_p(S_p)={\epsilon}_p(K(q_p))\cdot ({\theta}_p,discr(K(q_p)))_p$ $={\epsilon}_p(K(q_p))\cdot {\left( \frac{{\theta}_p}{p} \right)}^{log_p \mid A_p \mid}$ for odd $p$, and ${\epsilon}_2(S_2)=(-1)^{v_2+t_2(t_2-1)/2} \cdot  {\epsilon}_2(K_{\alpha}(q_2))\cdot ((-1)^{t_2}\cdot (-3)^{v_2},discr(K_{\alpha}(q_2)))_2$ $ =(-1)^{v_2+t_2(t_2-1)/2} \cdot  {\epsilon}_2(K_{\alpha}(q_2))\cdot (-1)^{\Lambda}$ $={\epsilon}_2(K_{\alpha}(q_2)) \cdot (-1)^{v_2\cdot (1+log_2 \mid A_2 \mid)}\cdot {\epsilon}_2$. Here we used $\Lambda = {log_2 (\mid A_2 \mid)}\cdot (3^{2v_2}-1)/8+(((-1)^{t_2}(-3)^{v_2}-1)/2)\cdot((\frac{discr(K_{\alpha}(q_2))}{\mid A_2 \mid}-1)/2)$ in order to simplify notation.\\ 

These conditions give the conditions 1-4 in the theorem.\\

Condition 5 in the theorem is equivalent to the condition that for some prime $p$ the p-adic lattice $S_p$ has no roots. It follows from the observation that $U^{(2)}(1)=\begin{pmatrix} 0 & 1\\ 1& 0 \end{pmatrix}$, $V^{(2)}(1)=\begin{pmatrix} 2 & 1\\ 1& 2 \end{pmatrix}$ and $K_{1}^{(p)}(1)^2=\begin{pmatrix} 1 & 0\\ 0& 1 \end{pmatrix}$ have roots. {\it QED}\\

{\bf Remark.} Theorem 9.1 gives a sufficient condition for the data $(t_{+},t_{-}, A_p,q_p)$ to correspond to an even lattice without roots. An example of the form $4x^2+4y^2+6z^2+8w^2$ shows that there exist even positive-definite forms without roots, which do have roots everywhere 'locally'. Indeed, according to \cite{Bhargava}, $1$ is the only positive integer, which is {\it not} represented by the integral quadratic form $2x^2+2y^2+3z^2+4w^2$. In particular, there is a cubic fourfold $X$ with $A_0(X)\cong
\begin{pmatrix}
4 & 0 & 0 & 0\\
0& 4 &0& 0\\
0& 0 &6& 0\\
0& 0 &0& 8
\end{pmatrix}$, even though all the p-adic localizations of this lattice have roots.\\

Let us check directly that for any $k\geq 1$ and any $N\in \mathbb Z$ congruence $2x^2+2y^2+3z^2+4w^2 \equiv N \; (mod\; p^k)$ has a solution. In fact, without loss of generality we may assume that $N>1$, because otherwise we can consider $N+M\cdot p^k$ instead of $N$ for sufficiently large $M$. Then according to Theorem 38 in \cite{Watson} (or by Hensel's lemma) it is sufficient to check that $2x^2+2y^2+3z^2+4w^2 \equiv N \; (mod\; 3)$ and $2x^2+2y^2+3z^2+4w^2 \equiv N \; (mod\; 16)$ are soluble. The first congruence has solution $x=0$, $y=0$, $z=0$, $w=1$ or $x=0$, $y=0$, $z=1$, $w=0$, if $N$ is a square $(mod\; 3)$, and solution $x=1$, $y=0$, $z=0$, $w=0$ otherwise. The existence of the solution of the second congruence can be checked by considering each of the $16$ possible values of $N \; (mod\; 16)$ separately. For the most interesting for us case ($N=1$), a solution is $x=2$, $y=1$, $z=1$, $w=1$.\\

This implies that $2x^2+2y^2+3z^2+4w^2 \equiv N \; (mod\; p^k)$ always has a solution. Hence for any $N\in \mathbb Z$ equation $2x^2+2y^2+3z^2+4w^2=N$ has a solution in p-adic integers. In particular, $2x^2+2y^2+3z^2+4w^2=1$, and so $4x^2+4y^2+6z^2+8w^2=2$ have solutions over ${\mathbb Z}_p$ for any $p$. This means that the lattice $A_0= \begin{pmatrix}
4 & 0 & 0 & 0\\
0& 4 &0& 0\\
0& 0 &6& 0\\
0& 0 &0& 8
\end{pmatrix}$ represents $2$ (i.e. has roots) everywhere locally. At the same time, it is clear that it has no roots (over $\mathbb Z$). Since it has signature $(\rho-1,0)$ with $\rho=5$, Corollary 5.1 implies that there exists a cubic fourfold $X$ such that the lattices $A_0$ and $A_0(X)$ are isomorphic.\\ 

\bibliographystyle{ams-plain}

\bibliography{IntersectionLattice}

\end{document}